\documentclass[12pt,twoside]{amsart}
\usepackage[latin1]{inputenc}
\usepackage{amsmath, amsthm, amscd, amsfonts, amssymb, graphicx}
\usepackage[bookmarksnumbered, plainpages]{hyperref}

\newtheorem{thm}{Theorem}[section]

\newtheorem{lem}[thm]{Lemma}

\newtheorem{defn}[thm]{Definition}

\numberwithin{equation}{section}

\allowdisplaybreaks
\begin{document}

\title{\bf A Kastler-Kalau-Walze type theorem for the $J$-twist of the Dirac operator with torsion}
\author{Jin Hong \hskip 0.4 true cm Siyao Liu \hskip 0.4 true cm  Yong Wang$^{*}$}

\thanks{{\scriptsize
\hskip -0.4 true cm \textit{2010 Mathematics Subject Classification:}
53C40; 53C42.
\newline \textit{$^{*}$Corresponding author}}}

\maketitle

\begin{abstract}
 \indent In this paper, we give a Lichnerowicz type formula for the $J$-twist of the Dirac operator with torsion. And we prove a Kastler-Kalau-Walze type theorem for the $J$-twist of the Dirac operator with torsion on 4-dimensional and 6-dimensional almost product Riemannian spin manifold with boundary.
\end{abstract}

\vskip 0.2 true cm


\pagestyle{myheadings}
\markboth{\rightline {\scriptsize Hong}}
         {\leftline{\scriptsize  A Kastler-Kalau-Walze type theorem for the $J$-twist of the Dirac operator with torsion}}

\bigskip
\bigskip


\section{ Introduction}
Noncommutative residues were first proposed in \cite{Gu,Wo}. Noncommutative residues are an important tool for studying noncommutative geometry. Therefore, much attention has been paid to noncommutative residues.
Connes put forward that the noncommutative residue of the square of the inverse of the Dirac operator was proportioned to the Einstein-Hilbert action, which is called the Kastler-Kalau-Walze theorem now\cite{Co1,Co2}. Kastler gave a brute-force proof of this theorem\cite{Ka}. In the same time, Kalau-Walze and Ackermann proved this theorem by using the normal coordinates system and the heat kernel expansion, respectively\cite{KW,Ac1}. The result of Connes was extended to the higher dimensional case \cite{U}. Wang generalized the Connes' results to the case of manifolds with boundary in \cite{Wa1,Wa2}, and proved the Kastler-Kalau-Walze type theorem for the Dirac operator and the signature operator on lower-dimensional manifolds with boundary. In \cite{Wa3}, for the Dirac operator, Wang computed $\widetilde{{\rm Wres}}[\pi^+D^{-1}\circ\pi^+D^{-1}]$, in these cases the boundary term vanished.
Most of the operators studied in the literature have the leading symbol $\sqrt{-1}c(\xi)$\cite{Wa4,Wa5,WW,WWW,Wa}. 

Earlier some preliminaries and lemmas about the Dirac operator $D$ and the $J$-twist are given in \cite{K}. In \cite{Chen1,Chen2}, the author got estimates on the higher eigenvalues of the Dirac operator on locally reducible Riemannian manifolds by the $J$-twist of the Dirac operator. It can be obtained by simple calculations that the leading symbol of the $J$-twist $D_{J}$ of the Dirac operator is not $\sqrt{-1}c(\xi)$.
Liu and Wang proved a Kastler-Kalau-Walze type theorem
for the J-twist $D_{J}$ of the Dirac operator on 3-dimensional, 4-dimensional and 6-dimensional almost
product Riemannian spin manifold with boundary in \cite{LW1,LW2}.
On the other hand, Jean-Michel Bismut \cite{BJM} proved a local index theorem for Dirac operators on a Riemannian manifold $M$ associated with connections on $TM$ which have non zero torsion.
In \cite{Ac2}, Ackermann and
Tolksdorf proved a generalized version of the well-known Lichnerowicz formula for the square of the most
general Dirac operator with torsion ${D_T}$ on an even-dimensional spin manifold associated to a metric connection with torsion.  Wang etc.\cite{WWJ} computed the lower dimensional volume $\widetilde{{\rm Wres}}[\pi^+{({D}_{T}^{*})}^{-p_1}\circ\pi^+{({D}_{T})}^{-p_2}]$ and got a Kastler-Kalau-Walze type theorems associated with Dirac operators
with torsion on compact manifolds with boundary.
So it is natural to consider the $J$-twist of the Dirac operator with torsion which is the common generalization of J-twisted Dirac operator and the Dirac operator with torsion. And it is natural to give the Kastler-Kalau-Walze type theorem for the J-twisted Dirac operator with torsion. The motivation of this paper is to give the Kastler-Kalau-Walze type theorem for the J-twisted Dirac operator with torsion.

The paper is structured as follows.
In Section 2, this paper will firstly introduce the basic notions of the almost product Riemannian manifold and the $J$-twist of the Dirac operator with torsion. We also give a Lichnerowicz type formula for the $J$-twist of the Dirac operator with torsion and a Kastler-Kalau-Walze type theorem for the $J$-twist of the Dirac operator with torsion on $n$-dimensional almost product Riemannian spin manifold without boundary. Next, we prove that the Kastler-Kalau-Walze type theorem on $4$-dimensional and $6$-dimensional almost product Riemannian spin manifold with boundary for the $J$-twist of the Dirac operator with torsion.

\vskip 1 true cm

\section{ The $J$-twist of the Dirac operator with torsion and its Lichnerowicz formula}

We give some definitions and basic notions which we will use in this paper.

Let $M$ be a smooth compact oriented Riemannian $n$-dimensional manifolds without boundary and $N$ be a vector bundle on $M$.
We say that $P$ is a differential operator of Laplace type, if it has locally the form
\begin{equation}\label{P}
P=-(g^{ij}\partial_i\partial_j+A^i\partial_i+B),
\end{equation}
where $\partial_{i}$  is a natural local frame on $TM,$ $(g^{ij})_{1\leq i,j\leq n}$ is the inverse matrix associated to the metric
matrix  $(g_{ij})_{1\leq i,j\leq n}$ on $M,$ $A^{i}$ and $B$ are smooth sections of $\textrm{End}(N)$ on $M$ (endomorphism).
If $P$ satisfies the form \eqref{P}, then there is a unique
connection $\nabla$ on $N$ and a unique endomorphism $E$ such that
 \begin{equation}
P=-[g^{ij}(\nabla_{\partial_{i}}\nabla_{\partial_{j}}- \nabla_{\nabla^{L}_{\partial_{i}}\partial_{j}})+E],\nonumber
\end{equation}
where $\nabla^{L}$ is the Levi-Civita connection on $M$. Moreover
(with local frames of $T^{*}M$ and $N$), $\nabla_{\partial_{i}}=\partial_{i}+\omega_{i} $
and $E$ are related to $g^{ij}$, $A^{i}$ and $B$ through
 \begin{eqnarray}
&&\omega_{i}=\frac{1}{2}g_{ij}\big(A^{i}+g^{kl}\Gamma_{ kl}^{j} \texttt{id}\big),\nonumber\\
&&E=B-g^{ij}\big(\partial_{i}(\omega_{j})+\omega_{i}\omega_{j}-\omega_{k}\Gamma_{ ij}^{k} \big),\nonumber
\end{eqnarray}
where $\Gamma_{ kl}^{j}$ is the  Christoffel coefficient of $\nabla^{L}$.

Let $M$ be a $n$-dimensional ($n\geq 3$) oriented compact Riemannian manifold with a Riemannian metric $g^{M}$.
We recall that the Dirac operator with torsion $D_T$ is locally given as following:
\begin{equation}
D_T=\sum_{i, j=1}^{n}g^{ij}c(\partial_{i})\nabla_{\partial_{j}}^{S(TM),T}=\sum_{i=1}^{n}c(e_{i})\nabla_{e_{i}}^{S(TM),T},\nonumber
\end{equation}
where $c(e_{i})$ be the Clifford action which satisfies the relation
\begin{align}
&c(e_{i})c(e_{j})+c(e_{j})c(e_{i})=-2g^{M}(e_{i}, e_{j})=-2\delta_i^j,\nonumber
\end{align}
\begin{align}
&\nabla_{e_{j}}^{S(TM),T}=e_{j}+\frac{1}{4}\sum_{l, t=1}^{n}\langle \nabla_{e_{j}}^{T}e_{l}, e_{t}\rangle c(e_{l})c(e_{t}),\nonumber
\end{align}
\begin{align}
&\langle \nabla_{e_{j}}^{T}e_{l}, e_{t}\rangle=\langle \nabla_{e_{j}}^{L}e_{l}, e_{t}\rangle + T(e_j, e_l, e_t)\nonumber
\end{align}
and
\begin{align}
	&\sigma^{j}=\frac{1}{4}\sum_{l, t=1}^{n}\langle \nabla_{e_{j}}^{L}e_{l}, e_{t}\rangle c(e_{l})c(e_{t}),\nonumber
\end{align}
where $T$ is a there form.

Let $J$ be a $(1, 1)$-tensor field on $(M, g^M)$ such that $J^2=\texttt{id},$
\begin{align}
&g^M(J(X), J(Y))=g^M(X, Y),\nonumber
\end{align}
for all vector fields $X,Y\in \Gamma(TM).$ Here $\texttt{id}$ stands for the identity map. $(M, g^M, J)$ is an almost product Riemannian manifold. We can also define on almost product Riemannian spin manifold the following $J$-twist of the Dirac operator with torsion $D_{J,T}$ by
\begin{align}
	&D_{J,T}:=\sum_{j=1}^{n}c[e_{j}]\nabla^{S(TM),T}_{J(e_{j})}=\sum_{j=1}^{n}c[J(e_{j})]\nabla^{S(TM),T}_{e_{j}} .
\end{align}

\begin{align}
	&T_{J}:=\frac{1}{4}\sum_{j, l, t =1}^{n}{T}\big(e_{j},e_{l}, e_{t}\big) c[J(e_{j})]c(e_{l})c(e_{t}).
\end{align}

Let $\xi=\sum_{j}\xi_{j}dx_{j},$  $\nabla^L_{\partial_{i}}\partial_{j}=\sum_{k}\Gamma_{ij}^{k}\partial_{k},$  we denote that
\begin{align}
 \Gamma^{k}=\sum_{i,j}g^{ij}\Gamma_{ij}^{k};\ \sigma^{j}=\sum_{i}g^{ij}\sigma_{i};\ \partial^{j}=\sum_{i}g^{ij}\partial_{i}.\nonumber
\end{align}

By (2.4), (2.5) in \cite{K}
\begin{align}
	D_{J, T}^{2}= & -\frac{1}{8} \sum_{i, j, k, l=1}^{n} R\left(J\left(e_{i}\right), J\left(e_{j}\right), e_{k}, e_{l}\right) c\left(e_{i}\right) c\left(e_{j}\right) c\left(e_{k}\right) c\left(e_{l}\right)-g^{i j} \partial_{i} \partial_{j}-2 \sigma^{j} \partial_{j} \\
	& +\Gamma^{k} \partial_{k} -g^{i j}\left[\partial_{i}\left(\sigma_{j}\right)+\sigma_{i} \sigma_{j}-\Gamma_{i j}^{k} \sigma_{k}\right]+\frac{1}{4} s+\sum_{\alpha, \beta=1}^{n} c\left[J\left(e_{\alpha}\right)\right] c\left[\left(\nabla_{e_{\alpha}}^{L} J\right) e_{\beta}\right] \nabla_{e_{\beta}}^{S} \nonumber\\
	& +\sum_{i, j}\left\{c\left[J\left(e_{i}\right)\right]\left\langle e_{i}, d x_{l}\right\rangle T_{J} \partial_{l}+T_{J}c\left[J\left(e_{i}\right)\right]\left\langle e_{i}, d x_{l}\right\rangle \partial_{l}\right\} \nonumber\\
	& +\sum_{i, j}\left\{c\left[J\left(e_{i}\right)\right]\left\langle e_{i}, d x_{l}\right\rangle \partial_{l}\left(T_{J}\right)+c\left[J\left(e_{i}\right)\right] \sigma_{i} T_{J}+T_{J} c\left[J\left(e_{i}\right)\right] \sigma_{i}\right\}+T_{J}^{2}.\nonumber
\end{align}
where $s$ is the scalar curvature, $e_i^*=g^{M}(e_i,\cdot)$ and $\langle X, dx_{k}\rangle=g^{M}(X, (dx_{k})^{*}),$ for a vector field $X.$
Then we have
\begin{align}
(\omega_{i})_{D_{J,T}^{2}}=&\sigma_{i}-\frac{1}{2}\sum_{\alpha,p=1}^{n}g_{ip}c[J(e_{\alpha})]c[(\nabla^{L}_{e_{\alpha}}J)(dx_{p})^{*}]\\
&-\frac{1}{2}\sum_{\alpha, p}g_{ip}c[J(e_{\alpha})]\left\langle e_{\alpha}, d x_{p}\right\rangle T_{J}-\frac{1}{2}\sum_{\alpha, p}g_{ip}T_{J}c[J(e_{\alpha})]\left\langle e_{\alpha}, d x_{p}\right\rangle\nonumber,
\end{align}
\begin{align}
E_{D_{J,T}^{2}}&=\frac{1}{8}\sum_{i,j,k,l=1}^{n}R(J(e_{i}), J(e_{j}), e_{k}, e_{l})c(e_{i})c(e_{j})c(e_{k})c(e_{l})+\sum_{i,j=1}^{n}g^{ij}[\partial_{i}(\sigma_{j})+\sigma_{i}\sigma_{j}-\Gamma^{k}_{ij}\sigma_{k}]\\
&-\sum_{\alpha,j=1}^{n}c[J(e_{\alpha})]c[(\nabla^{L}_{e_{\alpha}}J)(dx_{j})^{*}]\sigma_{j}-\sum_{i,j=1}^{n}g^{ij}\biggl\{\partial_{i} \bigg(\sigma_{j}-\frac{1}{2}\sum_{\nu,l=1}^{n}g_{jl}c[J(e_{\nu})]c[(\nabla^{L}_{e_{\nu}}J)(dx_{l})^{*}]\nonumber\\
&-\frac{1}{2}\sum_{v, l}^{n}g_{il}c[J(e_{v})]\left\langle e_{v}, d x_{l}\right\rangle T_{J}-\frac{1}{2}\sum_{v, l}^{n}g_{il}T_{J}c[J(e_{v})]\left\langle e_{v}, d x_{l}\right\rangle \bigg)\nonumber\\
&+\bigg(\sigma_{i}-\frac{1}{2}\sum_{\alpha,p=1}^{n}g_{ip}c[J(e_{\alpha})]c[(\nabla^{L}_{e_{\alpha}}J)(dx_{p})^{*}]-\frac{1}{2}\sum_{\alpha, p=1}^{n}g_{ip}c[J(e_{\alpha})]\left\langle e_{\alpha}, d x_{p}\right\rangle T_{J}\nonumber\\ 
&-\frac{1}{2}\sum_{\alpha, p=1}^{n}g_{ip}T_{J}c[J(e_{\alpha})]\left\langle e_{\alpha}, d x_{p}\right\rangle \bigg)\times \bigg(\sigma_{j}-\frac{1}{2}\sum_{\nu,l=1}^{n}g_{jl}c[J(e_{\nu})]c[(\nabla^{L}_{e_{\nu}}J)(dx_{l})^{*}]\nonumber\\
&-\frac{1}{2}\sum_{v, l=1}^{n}g_{jl}c[J(e_{v})]\left\langle e_{v}, d x_{l}\right\rangle T_{J}-\frac{1}{2}\sum_{v, l}^{n}g_{jl}T_{J}c[J(e_{v})]\left\langle e_{v}, d x_{l}\right\rangle \bigg)\nonumber\\
&-\sum_{k=1}^{n}\bigg(\sigma_{k}-\frac{1}{2}\sum_{\mu,h=1}^{n}g_{kh}c[J(e_{\mu})]c[(\nabla^{L}_{e_{\mu}}J)(dx_{h})^{*}]-\frac{1}{2}\sum_{\mu,h=1}^{n}g_{kh}c[J(e_{\mu})]\left\langle e_{\mu}, d x_{h}\right\rangle T_{J}\nonumber\\
&-\frac{1}{2}\sum_{\mu,h=1}^{n}g_{kh}T_{J}c[J(e_{\mu})]\left\langle e_{\mu},dx_{h}\right\rangle \bigg)\Gamma_{ij}^{k}\biggr\}-\frac{1}{4}s-T_{J}^{2}-\sum_{i=1}^{n}c[J(e_{i})]\sigma_{i} T_{J}-\sum_{i=1}^{n}T_{J}c[J(e_{i})]\sigma_{i}\nonumber \\
&-\sum_{\alpha, j}^{n}c[J(e_{\alpha})]\left\langle e_{\alpha}, d x_{j}\right\rangle \partial_{j}(T_{J}) .\nonumber
\end{align}

Since $E$ is globally defined on $M$, taking normal coordinates at $x_0$, we have $\sigma^{i}(x_0)=0$, $\Gamma^k(x_0)=0,$ $g^{ij}(x_0)=\delta^j_i,$ $\partial^{j}(x_0)=e_{j},$ $\partial^{j}[c(\partial_{j})](x_0)=0,$ $\nabla^{L}_{e_{j}}e_{k}(x_0)=0$ and $\nabla^{S}_{Y}(c(X))=c(\nabla^{L}_{Y}X),$ a simple calculation shows that
\begin{align}
&\frac{1}{2}\sum_{\nu,j=1}^{n}\partial^{j}(c[J(e_{\nu})]c[(\nabla^{L}_{e_{\nu}}J)(dx_{j})^{*}])(x_0)\nonumber\\
&=\frac{1}{2}\sum_{\nu,j=1}^{n}c[\nabla_{e_{j}}^{L}(J)e_{\nu}]c[(\nabla^{L}_{e_{\nu}}J)e_{j}](x_0)+\frac{1}{2}\sum_{\nu,j=1}^{n}c[J(e_{\nu})]c[(\nabla^{L}_{e_{j}}(\nabla^{L}_{e_{\nu}}(J)))e_{j}-(\nabla^{L}_{\nabla^{L}_{e_{j}}e_{\nu}}(J))e_{j}](x_0),\nonumber
\end{align}

so that
\begin{align}
&E_{D_{J,T}^{2}}(x_0)=\frac{1}{8}\sum_{i,j,k,l=1}^{n}R(J(e_{i}), J(e_{j}), e_{k}, e_{l})c(e_{i})c(e_{j})c(e_{k})c(e_{l})\nonumber\\
&+\frac{1}{2}\sum_{\nu,j=1}^{n}c[\nabla_{e_{j}}^{L}(J)e_{\nu}]c[(\nabla^{L}_{e_{\nu}}J)e_{j}]+\frac{1}{2}\sum_{\nu,j=1}^{n}c[J(e_{\nu})]c[(\nabla^{L}_{e_{j}}(\nabla^{L}_{e_{\nu}}(J)))e_{j}-(\nabla^{L}_{\nabla^{L}_{e_{j}}e_{\nu}}(J))e_{j}]\nonumber\\
&-\frac{1}{4}\sum_{\alpha,\nu,j=1}^{n}c[J(e_{\alpha})]c[(\nabla^{L}_{e_{\alpha}}J)e_{j}]c[J(e_{\nu})]c[(\nabla^{L}_{e_{\nu}}J)e_{j}]-\frac{1}{4}s-T_{J}^{2}\nonumber\\
&-\frac{1}{4} \sum_{\hat{j},\hat{l},\hat{\alpha},j=1}^{n}\left[J\left(e_{j}\right)\right]T\left(e_{\hat{j}},e_{\hat{l}}, e_{\hat{\alpha}}\right) c\left[\nabla_{\partial _{j}}^{L}(J) e_{\hat{j}}\right] c\left(e_{\hat{l}}\right)c\left(e_{\hat{\alpha}}\right)\nonumber\\
&+\frac{1}{4}\sum_{\hat{j},\hat{l},\hat{\alpha},j=1}^{n}T\left(e_{\hat{j}},e_{\hat{l}}, e_{\hat{\alpha}}\right) \biggl\{ c\left[\nabla_{e _{j}}^{L}(J) e_{j}\right]c[J(e_{\hat{j}})]c(e_{\hat{l}})c(e_{\hat{\alpha}})+c[J(e_{j})]c\left[\nabla_{e _{j}}^{L}(J) e_{\hat{j}}\right]c(e_{\hat{l}})c(e_{\hat{\alpha}}) \biggr\}\nonumber\\
&-\frac{1}{4} \sum_{\alpha,j=1}^{n}c[J(e_{\alpha})]c\left[\nabla_{e _{\alpha}}^{L}(J) e_{j}\right]\biggl\{ c[J(dx_{j})]T_J+T_Jc[J(dx_{j})]\biggr\} \nonumber\\
&-\frac{1}{4} \sum_{v,\alpha,j=1}^{n}c[J(dx_{j})T_J \biggl\{c[J(e_{v})c\left[\nabla_{e _{v}}^{L}(J) e_{j}\right]+c[J(dx_{j})T_J+T_Jc[J(dx_{j})] \biggr\}\nonumber\\
&-\frac{1}{4} \sum_{v,\alpha,j=1}^{n} T_J c[J(dx_{j}) \biggl\{c[J(e_{v})c\left[\nabla_{e _{v}}^{L}(J) e_{j}\right]+c[J(dx_{j})T_J+T_Jc[J(dx_{j})] \biggr\} .\nonumber
\end{align}

We get the following Lichnerowicz formulas:
\begin{align}
D_{J,T}^{2}
=&-g^{ij}(\nabla_{\partial_{i}}\nabla_{\partial_{j}}-\nabla_{\nabla^{L}_{\partial_{i}}\partial_{j}})-\frac{1}{8}\sum_{i,j,k,l=1}^{n}R(J(e_{i}), J(e_{j}), e_{k}, e_{l})c(e_{i})c(e_{j})c(e_{k})c(e_{l})+\frac{1}{4}s+T_{J}^{2}\\
&-\frac{1}{2}\sum_{\nu,j=1}^{n}c[\nabla_{e_{j}}^{L}(J)e_{\nu}]c[(\nabla^{L}_{e_{\nu}}J)e_{j}]-\frac{1}{2}\sum_{\nu,j=1}^{n}c[J(e_{\nu})]c[(\nabla^{L}_{e_{j}}(\nabla^{L}_{e_{\nu}}(J)))e_{j}-(\nabla^{L}_{\nabla^{L}_{e_{j}}e_{\nu}}(J))e_{j}]\nonumber\\
&+\frac{1}{4}\sum_{\alpha,\nu,j=1}^{n}c[J(e_{\alpha})]c[(\nabla^{L}_{e_{\alpha}}J)e_{j}]c[J(e_{\nu})]c[(\nabla^{L}_{e_{\nu}}J)e_{j}]\nonumber\\
&-\frac{1}{4} \sum_{\hat{j},\hat{l},\hat{\alpha},j=1}^{n}\left[J\left(e_{j}\right)\right] T\big( J(e_{\hat{j}}), e_{\hat{l}}, e_{\hat{\alpha}}\big) c\left[\nabla_{\partial _{j}}^{L}(J) e_{\hat{j}}\right] c\left(e_{\hat{l}}\right)c\left(e_{\hat{\alpha}}\right)\nonumber\\
&+\frac{1}{4}\sum_{\hat{j},\hat{l},\hat{\alpha},j=1}^{n} T\big( J(e_{\hat{j}}), e_{\hat{l}}, e_{\hat{\alpha}}\big)\biggl\{ c\left[\nabla_{e _{j}}^{L}(J) e_{j}\right]c[J(e_{\hat{j}})]c(e_{\hat{l}})c(e_{\hat{\alpha}})+c[J(e_{j})]c\left[\nabla_{e _{j}}^{L}(J) e_{\hat{j}}\right]c(e_{\hat{l}})c(e_{\hat{\alpha}}) \biggr\}\nonumber\\
&-\frac{1}{4} \sum_{\alpha,j=1}^{n}c[J(e_{\alpha})]c\left[\nabla_{e _{\alpha}}^{L}(J) e_{j}\right]\biggl\{ c[J(dx_{j})]T_J+T_Jc[J(dx_{j})]\biggr\} \nonumber\\
&-\frac{1}{4} \sum_{v,\alpha,j=1}^{n}c[J(dx_{j})T_J \biggl\{c[J(e_{v})c\left[\nabla_{e _{v}}^{L}(J) e_{j}\right]+c[J(dx_{j})T_J+T_Jc[J(dx_{j})] \biggr\}\nonumber\\
&-\frac{1}{4} \sum_{v,\alpha,j=1}^{n} T_J c[J(dx_{j}) \biggl\{c[J(e_{v})c\left[\nabla_{e _{v}}^{L}(J) e_{j}\right]+c[J(dx_{j})T_J+T_Jc[J(dx_{j})] \biggr\} .\nonumber
\end{align}

Now, let's calculate the specification of $D_{J,T}^{3}$.
\begin{align}\label{D3}
	D_{J,T}^{3} & =\sum_{i=1}^{n} c\left[J\left(e_{i}\right)\right]\left\langle e_{i}, d x_{l}\right\rangle\left(-g^{i j} \partial_{l} \partial_{i} \partial_{j}\right)+\sum_{i=1}^{n} c\left[J\left(e_{i}\right)\right]\left\langle e_{i}, d x_{l}\right\rangle\bigg[-\left(\partial_{l} g^{i j}\right) \partial_{i} \partial_{j}-g^{i j}\left(2 \sigma_{i} \partial_{j}\right.\\
	& \left. - \Gamma_{i j}^{k} \partial_{k}\right) \partial_{l}+ \sum_{\alpha, \beta, k=1}^{n} c\left[J\left(e_{\alpha}\right)\right] c\left[\left(\nabla_{e_{\alpha}}^{L} J\right) e_{\beta}\right]\left\langle e_{\beta}, d x_{k}\right\rangle \partial_{l} \partial_{k}\bigg]+\sum_{i=1}^{n} c\left[J\left(e_{i}\right)\right] \sigma_{i}\left(-g^{i j} \partial_{i} \partial_{j}\right)\nonumber \\
	& -T_{J} g^{i j} \partial_{i} \partial_{j}+ \sum_{i=1}^{n} c\left[J\left(e_{i}\right)\right]\left\langle e_{i}, d x_{l}\right\rangle\left(c\left[J\left(e_{i}\right)\right] T_{J}+T_{J} c\left[J\left(e_{i}\right)\right]\right) \partial_{i} \partial_{j} \nonumber\\
	& +\sum_{i=1}^{n} c\left[J\left(e_{i}\right)\right]\left\langle e_{i}, d x_{l}\right\rangle\bigg[-2\left(\partial_{l} g^{i j}\right) \sigma_{i} \partial_{j}-2 g^{i j}\left(\partial_{l} \sigma_{i}\right) \partial_{j}+\left(\partial_{l} g^{i j}\right) \Gamma_{i j}^{k} \partial_{k}+g^{i j}\left(\partial_{l} \Gamma_{i j}^{k}\right) \partial_{k} \nonumber\\
	& +\sum_{\alpha, \beta, k=1}^{n} \partial_{l}\left(c\left[J\left(e_{\alpha}\right)\right] c\left[\left(\nabla_{e_{\alpha}}^{L} J\right) e_{\beta}\right]\right)\left\langle e_{\beta}, d x_{k}\right\rangle \partial_{k}+\sum_{\alpha, \beta, k=1}^{n} c\left[J\left(e_{\alpha}\right)\right] c\left[\left(\nabla_{e_{\alpha}}^{L} J\right) e_{\beta}\right]\left(\partial_{l}\left\langle e_{\beta}, d x_{k}\right\rangle\right) \partial_{k} \nonumber\\
	& +\left(\partial_{l} c\left[J\left(e_{i}\right)\right]\right) T_{J}+c\left[J\left(e_{i}\right)\right]\left(\partial_{l} T_{J}\right)+\left(\partial_{l} T_{J}\right) c\left[J\left(e_{i}\right)\right]+T_{J}\left(\partial_{l} c\left[J\left(e_{i}\right)\right]\right)\bigg] \nonumber\\
	& +\sum_{i=1}^{n} c\left[J\left(e_{i}\right)\right]\left\langle e_{i}, d x_{l}\right\rangle \partial_{l}\bigg[-\frac{1}{8} \sum_{i, j, k, l=1}^{n} R\left(J\left(e_{i}\right), J\left(e_{j}\right), e_{k}, e_{l}\right) c\left(e_{i}\right) c\left(e_{j}\right) c\left(e_{k}\right) c\left(e_{l}\right)-g^{i j}\left(\left(\partial_{i} \sigma_{j}\right)\right.\nonumber \\
	& \left.+\sigma_{i} \sigma_{j}-\Gamma_{i j}^{k} \sigma_{k}\right)+\frac{1}{4} s+\sum_{\alpha, \beta, k=1}^{n} c\left[J\left(e_{\alpha}\right)\right] c\left[\left(\nabla_{e_{\alpha}}^{L} J\right) e_{\beta}\right]\left\langle e_{\beta}, d x_{k}\right\rangle \sigma_{k}+c\left[J\left(e_{j}\right)\right]\left(\partial_{j} T_{J}\right)\nonumber\\
	&+c\left[J\left(e_{\alpha}\right)\right] \sigma_{\alpha} T_{J}  +T_{J} c\left[J\left(e_{\alpha}\right)\right] \sigma_{\alpha}+T_{J}^{2}\bigg]+\left[\sum_{i=1}^{n} c\left[J\left(e_{i}\right)\right] \sigma_{i}+T_{J}\right]\nonumber\\
	&\times \bigg[\frac{1}{4} s+\left(c\left[J\left(e_{j}\right)\right] T_{J}+T_{J} c\left[J\left(e_{j}\right)\right]\right) \partial_{j} -\frac{1}{8} \sum_{i, j, k, l=1}^{n} R\left(J\left(e_{i}\right), J\left(e_{j}\right), e_{k}, e_{l}\right) c\left(e_{i}\right) c\left(e_{j}\right) c\left(e_{k}\right) c\left(e_{l}\right)\nonumber \\
	&-2 \sigma^{j} \partial_{j}+\Gamma^{k} \partial_{k}-g^{i j}\left(\left(\partial_{i} \sigma_{j}\right)+\sigma_{i} \sigma_{j}\right. \left.-\Gamma_{i j}^{k} \sigma_{k}\right)+\sum_{\alpha, \beta, k=1}^{n} c\left[J\left(e_{\alpha}\right)\right] c\left[\left(\nabla_{e_{\alpha}}^{L} J\right) e_{\beta}\right]\left\langle e_{\beta}, d x_{k}\right\rangle\left(\partial_{k}+\sigma_{k}\right) \nonumber\\
	& +c\left[J\left(e_{j}\right)\right]\left(\partial_{j} T_{J}\right)+c\left[J\left(e_{\alpha}\right)\right] \sigma_{\alpha} T_{J}+T_{J} c\left[J\left(e_{\alpha}\right)\right] \sigma_{\alpha}+T_{J}^{2}\bigg] .\nonumber
\end{align}

According to the detailed descriptions in \cite{Ac1}, we know that the noncommutative residue of a generalized laplacian $\widetilde{\Delta}$ is expressed as
\begin{equation}
(n-2)\Phi_{2}(\widetilde{\Delta})=\pi^{-\frac{n}{2}}\Gamma(\frac{n}{2})\widetilde{res}(\widetilde{\Delta}^{-\frac{n}{2}+1}),\nonumber
\end{equation}
where $\Phi_{2}(\widetilde{\Delta})$ denotes the integral over the diagonal part of the second
coefficient of the heat kernel expansion of $\widetilde{\Delta}$.
Now let $\widetilde{\Delta}=D_{J,T}^{2}$. Since $D_{J,T}^{2}$ is a generalized laplacian, we can suppose ${D_{J}}^{2}=\Delta-E$, then, we have
\begin{align}\label{r1}
{\rm Wres}(D_{J,T}^{2})^{-\frac{n-2}{2}}
=\frac{(n-2)\pi^{\frac{n}{2}}}{(\frac{n}{2}-1)!}\int_{M}{\rm tr}(-\frac{1}{6}s+E_{D_{J,T}^{2}})d{\rm Vol_{M} },
\end{align}
where ${\rm Wres}$ denote the noncommutative residue, ${\rm tr}$ denote ${\rm tr}$. By
\begin{align}\label{r2}
{\rm tr}[c(X)c(Y)]=-g^{M}(X, Y){\rm tr}[\texttt{id}];
\end{align}
\begin{align}\label{r3}
{\rm tr}[c(X)c(Y)c(Z)c(W)]&=g^{M}(X, W)g^{M}(Y, Z){\rm tr}[\texttt{id}]-g^{M}(X, Z)g^{M}(Y, W){\rm tr}[\texttt{id}]\\
&+g^{M}(X, Y)g^{M}(Z, W){\rm tr}[\texttt{id}];\nonumber
\end{align}
\begin{align}\label{r4}
	&{\rm tr}[c(X)c(Y)c(Z)c(W)c(V)c(U)]\\
	=&-g^{M}(X, U)g^{M}(Y, V)g^{M}(Z, W){\rm tr}[\texttt{id}]+g^{M}(X, U)g^{M}(Y, W)g^{M}(Z, V){\rm tr}[\texttt{id}]\nonumber\\
	&-g^{M}(X, U)g^{M}(Y, Z)g^{M}(W, V){\rm tr}[\texttt{id}]+g^{M}(X,V)g^{M}(Y, U)g^{M}(Z,W){\rm tr}[\texttt{id}]\nonumber\\
	&-g^{M}(X,V)g^{M}(Y, W)g^{M}(Z,U){\rm tr}[\texttt{id}]+g^{M}(X,V)g^{M}(Y,Z)g^{M}(W,U){\rm tr}[\texttt{id}]\nonumber\\
	&-g^{M}(X,W)g^{M}(Y,U)g^{M}(Z,V){\rm tr}[\texttt{id}]+g^{M}(X,W)g^{M}(Y,V)g^{M}(Z,U){\rm tr}[\texttt{id}]\nonumber\\
	&-g^{M}(X,W)g^{M}(Y,Z)g^{M}(V,U){\rm tr}[\texttt{id}]+g^{M}(X,Z)g^{M}(Y,U)g^{M}(W,V){\rm tr}[\texttt{id}]\nonumber\\
	&-g^{M}(X,Z)g^{M}(Y,V)g^{M}(W,U){\rm tr}[\texttt{id}]+g^{M}(X,Z)g^{M}(Y,W)g^{M}(V,U){\rm tr}[\texttt{id}]\nonumber\\
	&-g^{M}(X,Y)g^{M}(Z,U)g^{M}(W,V){\rm tr}[\texttt{id}]+g^{M}(X,Y)g^{M}(Z,V)g^{M}(W,U){\rm tr}[\texttt{id}]\nonumber\\
	&-g^{M}(X,Y)g^{M}(Z,W)g^{M}(V,U){\rm tr}[\texttt{id}]\nonumber
\end{align}
where $X,Y,Z,W,V,U\in \Gamma(TM).$ We calculate that
\begin{align}\label{r5}
&\sum_{i,j,k,l=1}^{n}{\rm tr}[R(J(e_{i}), J(e_{j}), e_{k}, e_{l})c(e_{i})c(e_{j})c(e_{k})c(e_{l})]\\
&=\sum_{i,j,k,l=1}^{n}R(J(e_{i}), J(e_{j}), e_{k}, e_{l}){\rm tr}[c(e_{i})c(e_{j})c(e_{k})c(e_{l})]\nonumber\\
&=\sum_{i,j=1}^{n}R(J(e_{i}), J(e_{j}), e_{j}, e_{i}){\rm tr}[\texttt{id}]-\sum_{i,j=1}^{n}R(J(e_{i}), J(e_{j}), e_{i}, e_{j}){\rm tr}[\texttt{id}]\nonumber\\
&=2\sum_{i,j=1}^{n}R(J(e_{i}), J(e_{j}), e_{j}, e_{i}){\rm tr}[\texttt{id}],\nonumber
\end{align}
\begin{align}\label{r6}
&\sum_{\nu,j=1}^{n}{\rm tr}[c[\nabla_{e_{j}}^{L}(J)e_{\nu}]c[(\nabla^{L}_{e_{\nu}}J)e_{j}]]=-\sum_{\nu,j=1}^{n}g^{M}(\nabla_{e_{j}}^{L}(J)e_{\nu}, (\nabla^{L}_{e_{\nu}}J)e_{j}){\rm tr}[\texttt{id}],
\end{align}
\begin{align}\label{r7}
&\sum_{\nu,j=1}^{n}{\rm tr}[c[J(e_{\nu})]c[(\nabla^{L}_{e_{j}}(\nabla^{L}_{e_{\nu}}(J)))e_{j}-(\nabla^{L}_{\nabla^{L}_{e_{j}}e_{\nu}}(J))e_{j}]]\\
&=-\sum_{\nu,j=1}^{n}g^{M}(J(e_{\nu}), (\nabla^{L}_{e_{j}}(\nabla^{L}_{e_{\nu}}(J)))e_{j}-(\nabla^{L}_{\nabla^{L}_{e_{j}}e_{\nu}}(J))e_{j}){\rm tr}[\texttt{id}],\nonumber
\end{align}
\begin{align}\label{r8}
&\sum_{\alpha,\nu,j=1}^{n}{\rm tr}[c[J(e_{\alpha})]c[(\nabla^{L}_{e_{\alpha}}J)e_{j}]c[J(e_{\nu})]c[(\nabla^{L}_{e_{\nu}}J)e_{j}]]\\
&=\sum_{\alpha,\nu,j=1}^{n}g^{M}(J(e_{\alpha}), (\nabla^{L}_{e_{\nu}}J)e_{j})g^{M}((\nabla^{L}_{e_{\alpha}}J)e_{j}, J(e_{\nu})){\rm tr}[\texttt{id}]\nonumber\\
&-\sum_{\nu,j=1}^{n}g^{M}((\nabla^{L}_{e_{\nu}}J)e_{j}, (\nabla^{L}_{e_{\nu}}J)e_{j}){\rm tr}[\texttt{id}]\nonumber\\
&+\sum_{\alpha,\nu,j=1}^{n}g^{M}(J(e_{\alpha}), (\nabla^{L}_{e_{\alpha}}J)e_{j})g^{M}(J(e_{\nu}), (\nabla^{L}_{e_{\nu}}J)e_{j}){\rm tr}[\texttt{id}].\nonumber
\end{align}

By applying the formulae shown in \eqref{r1}-\eqref{r8}, we get:
\begin{thm} If $M$ is a $n$-dimensional almost product Riemannian spin manifold without boundary, we have the following:
\begin{align}
&{\rm Wres}(D_{J,T}^{-n+2})=\nonumber\\
&\frac{(n-2)\pi^{\frac{n}{2}}}{(\frac{n}{2}-1)!}\int_{M}2^{\frac{n}{2}} \biggl\{ (\frac{1}{4}\sum_{i,j=1}^{n}R(J(e_{i}), J(e_{j}), e_{j}, e_{i})
-\frac{1}{2}\sum_{\nu,j=1}^{n}g^{M}(\nabla_{e_{j}}^{L}(J)e_{\nu}, (\nabla^{L}_{e_{\nu}}J)e_{j})\nonumber\\
&-\frac{1}{2}\sum_{\nu,j=1}^{n}g^{M}(J(e_{\nu}), (\nabla^{L}_{e_{j}}(\nabla^{L}_{e_{\nu}}(J)))e_{j}-(\nabla^{L}_{\nabla^{L}_{e_{j}}e_{\nu}}(J))e_{j})\nonumber\\
&-\frac{1}{4}\sum_{\alpha,\nu,j=1}^{n}g^{M}(J(e_{\alpha}), (\nabla^{L}_{e_{\nu}}J)e_{j})g^{M}((\nabla^{L}_{e_{\alpha}}J)e_{j}, J(e_{\nu}))\nonumber\\
&-\frac{1}{4}\sum_{\alpha,\nu,j=1}^{n}g^{M}(J(e_{\alpha}), (\nabla^{L}_{e_{\alpha}}J)e_{j})g^{M}(J(e_{\nu}), (\nabla^{L}_{e_{\nu}}J)e_{j})+\frac{1}{4}\sum_{\nu,j=1}^{n}g^{M}((\nabla^{L}_{e_{\nu}}J)e_{j}, (\nabla^{L}_{e_{\nu}}J)e_{j})-\frac{5}{12}s\nonumber\\
&-\frac{1}{2}\sum_{j, m, p} T\left(e_m, e_{p}, J\left(e_{m}\right)\right) g^{M}\left((\nabla_{e_{j}}^{L}J) e_{j} ,e_{p}\right)-\frac{1}{2}\sum_{j, m, \alpha} T\left(e_m, J\left(e_{j}\right), J\left(e_{\alpha}\right)\right) g^{M}\left((\nabla_{e_{\alpha}}^{L}J) e_{j} ,J(e_m)\right)\nonumber\\
&-\frac{1}{2}\sum_{j, m, \alpha} T\left(e_m, J\left(e_{j}\right), J\left(e_{m}\right)\right) g^{M}\left((\nabla_{e_{\alpha}}^{L}J) e_{j} ,J(e_{\alpha})\right)+\frac{1}{2}\sum_{m,p, \alpha} T\left(e_m, e_p, J\left(e_{\alpha}\right)\right) g^{M}\left((\nabla_{e_{\alpha}}^{L}J) e_{m} ,e_p\right)\nonumber\\
&-\frac{1}{2}\sum_{m,p, \alpha} T\left(e_m, e_p, J\left(e_{\alpha}\right)\right) g^{M}\left((\nabla_{e_{m}}^{L}J) e_{\alpha} ,e_p\right)-\frac{1}{4} \sum_{j,p,m} T(e_j, e_p, J(e_m))T(e_m, e_p, J(e_j))\nonumber\\
&-\frac{1}{2} \sum_{j, l, \alpha} T^2\left(J\left(e_{j}\right), e_{l}, e_{\alpha}\right)\biggr\}d{\rm Vol_M}\nonumber,
\end{align}
where $s$ is the scalar curvature.
\end{thm}

\section{ A Kastler-Kalau-Walze type theorem for $4$-dimensional manifolds with boundary}
\indent Firstly, we need the basic notions of  Boutet de Monvel's calculus and the definition of the noncommutative residue for manifolds with boundary that will be used throughout the paper. For the details, see Ref.\cite{Wa3}. And we use the same definitions and symbols as in \cite{Wa3}.\\
\indent Let $U\subset M$ be a collar neighborhood of $\partial M$ which is diffeomorphic with $\partial M\times [0,1)$. By the definition of $h(x_n)\in C^{\infty}([0,1))$
and $h(x_n)>0$, there exists $\widehat{h}\in C^{\infty}((-\varepsilon,1))$ such that $\widehat{h}|_{[0,1)}=h$ and $\widehat{h}>0$ for some
sufficiently small $\varepsilon>0$. Then there exists a metric $g'$ on $\widetilde{M}=M\bigcup_{\partial M}\partial M\times
(-\varepsilon,0]$ which has the form on $U\bigcup_{\partial M}\partial M\times (-\varepsilon,0 ]$
\begin{equation}
g'=\frac{1}{\widehat{h}(x_{n})}g^{\partial M}+dx _{n}^{2} ,\nonumber
\end{equation}
such that $g'|_{M}=g$. We fix a metric $g'$ on the $\widetilde{M}$ such that $g'|_{M}=g$.

\begin{defn}{\rm\cite{Wa3} }
Lower dimensional volumes of almost product spin manifolds with boundary are defined by
 \begin{equation}
{\rm Vol}^{(p_1,p_2)}_nM:= \widetilde{{\rm Wres}}[\pi^+D_{J,T}^{-p_1}\circ\pi^+D_{J,T}^{-p_2}],\nonumber
\end{equation}
where $\widetilde{{\rm Wres}}$ denote the noncommutative residue for manifolds with boundary as in \cite{FGLS}.
\end{defn}
We can get the spin structure on $\widetilde{M}$ by extending the spin structure on $M.$ Let $D$ be the Dirac operator associated to $g'$ on the spinors bundle $S(T\widetilde{M}).$ 

By \cite{Wa3}, we get
\begin{align}\label{w}
\widetilde{{\rm Wres}}[\pi^+D_{J,T}^{-p_1}\circ\pi^+D_{J,T}^{-p_2}]=\int_M\int_{|\xi|=1}{\rm
trace}_{S(TM)}[\sigma_{-n}(D_{J,T}^{-p_1-p_2})]\sigma(\xi)dx+\int_{\partial M}\Phi
\end{align}
and
\begin{align} \label{phi}
\Phi&=\int_{|\xi'|=1}\int^{+\infty}_{-\infty}\sum^{\infty}_{j, k=0}\sum\frac{(-i)^{|\alpha|+j+k+1}}{\alpha!(j+k+1)!}
\times {\rm tr}_{S(TM)}[\partial^j_{x_n}\partial^\alpha_{\xi'}\partial^k_{\xi_n}\sigma^+_{r}(D_{J,T}^{-p_1})(x',0,\xi',\xi_n)\\
&\times\partial^\alpha_{x'}\partial^{j+1}_{\xi_n}\partial^k_{x_n}\sigma_{l}(D_{J,T}^{-p_2})(x',0,\xi',\xi_n)]d\xi_n\sigma(\xi')dx',\nonumber
\end{align}
 where the sum is taken over $r+l-k-|\alpha|-j-1=-n,~~r\leq -p_1,l\leq -p_2$.
 
 Since $[\sigma_{-n}(D_{J,T}^{-p_1-p_2})]|_M$ has the same expression as $\sigma_{-n}(D_{J,T}^{-p_1-p_2})$ in the case of manifolds without
boundary, so locally we can compute the first term by \cite{Ka}, \cite{KW}, \cite{Wa3}, \cite{Po}.

We firstly compute
\begin{equation}\label{wt}
\widetilde{{\rm Wres}}[\pi^+{{D}_{J,T}^{-1}}\circ\pi^+{{D}_{J,T}^{-1}}]=\int_M\int_{|\xi|=1}{\rm
trace}_{S(TM)}[\sigma_{-4}({{D}_{J,T}^{-2}})]\sigma(\xi)dx+\int_{\partial M}\Psi,
\end{equation}
where
\begin{align}\label{phit}
\Psi &=\int_{|\xi'|=1}\int^{+\infty}_{-\infty}\sum^{\infty}_{j, k=0}\sum\frac{(-i)^{|\alpha|+j+k+1}}{\alpha!(j+k+1)!}
\times {\rm tr}_{S(TM)}[\partial^j_{x_n}\partial^\alpha_{\xi'}\partial^k_{\xi_n}\sigma^+_{r}({{D}_{J,T}^{-1}})\\
&(x',0,\xi',\xi_n)\times\partial^\alpha_{x'}\partial^{j+1}_{\xi_n}\partial^k_{x_n}\sigma_{l}({{D}_{J,T}^{-1}})(x',0,\xi',\xi_n)]d\xi_n\sigma(\xi')dx',\nonumber
\end{align}
the sum is taken over $r+l-k-j-|\alpha|-1=-4, r\leq -1, l\leq-1$.\\

\indent
Then we give the interior of $\widetilde{{\rm Wres}}[\pi^+{D}_{J,T}^{-1}\circ\pi^+{D}_{J,T}^{-1}]$,
\begin{align}
&\int_M\int_{|\xi|=1}{\rm tr}_{S(TM)}[\sigma_{-4}({D}_{J,T}^{-2})]\sigma(\xi)dx=2\pi^{2}\\
&\int_{M}\biggl\{ \sum_{i,j=1}^{n}R(J(e_{i}), J(e_{j}), e_{j}, e_{i})
-2\sum_{\nu,j=1}^{n}g^{M}(\nabla_{e_{j}}^{L}(J)e_{\nu}, (\nabla^{L}_{e_{\nu}}J)e_{j})\nonumber \\
&-2\sum_{\nu,j=1}^{n}g^{M}(J(e_{\nu}), (\nabla^{L}_{e_{j}}(\nabla^{L}_{e_{\nu}}(J)))e_{j}-(\nabla^{L}_{\nabla^{L}_{e_{j}}e_{\nu}}(J))e_{j})\nonumber\\
&-\sum_{\alpha,\nu,j=1}^{n}g^{M}(J(e_{\alpha}), (\nabla^{L}_{e_{\nu}}J)e_{j})g^{M}((\nabla^{L}_{e_{\alpha}}J)e_{j}, J(e_{\nu}))\nonumber\\
&-\sum_{\alpha,\nu,j=1}^{n}g^{M}(J(e_{\alpha}), (\nabla^{L}_{e_{\alpha}}J)e_{j})g^{M}(J(e_{\nu}), (\nabla^{L}_{e_{\nu}}J)e_{j})+\sum_{\nu,j=1}^{n}g^{M}((\nabla^{L}_{e_{\nu}}J)e_{j}, (\nabla^{L}_{e_{\nu}}J)e_{j})-\frac{5}{3}s\nonumber\\
&-2\sum_{j, m, p} T\left(e_m, e_{p}, J\left(e_{m}\right)\right) g^{M}\left((\nabla_{e_{j}}^{L}J) e_{j} ,e_{p}\right)-2\sum_{j, m, \alpha} T\left(e_m, J\left(e_{j}\right), J\left(e_{\alpha}\right)\right) g^{M}\left((\nabla_{e_{\alpha}}^{L}J) e_{j} ,J(e_m)\right)\nonumber\\
&-2\sum_{j, m, \alpha} T\left(e_m, J\left(e_{j}\right), J\left(e_{m}\right)\right) g^{M}\left((\nabla_{e_{\alpha}}^{L}J) e_{j} ,J(e_{\alpha})\right)+2\sum_{m,p, \alpha} T\left(e_m, e_p, J\left(e_{\alpha}\right)\right) g^{M}\left((\nabla_{e_{\alpha}}^{L}J) e_{m} ,e_p\right)\nonumber\\
&-2\sum_{m,p, \alpha} T\left(e_m, e_p, J\left(e_{\alpha}\right)\right) g^{M}\left((\nabla_{e_{m}}^{L}J) e_{\alpha} ,e_p\right)- \sum_{j,p,m} T(e_j, e_p, J(e_m))T(e_m, e_p, J(e_j))\nonumber\\
&-2 \sum_{j, l, \alpha} T^2\left(J\left(e_{j}\right), e_{l}, e_{\alpha}\right)\biggr\}d{\rm Vol_M}\nonumber
\end{align}
where $\Omega_{n}=\frac{2\pi^\frac{n}{2}}{\Gamma(\frac{n}{2})}.$

\indent We begin by computing $\int_{\partial M} \Psi$. Since, some operators have the following symbols.
\begin{lem} The following identities hold:
\begin{align}
\sigma_1({D}_{J,T})=&\sqrt{-1}c[J(\xi)];\nonumber\\
\sigma_0({D}_{J,T})=&
-\frac{1}{4}\sum_{i,j,k=1}^{n}\omega_{j,k}(e_i)c[J(e_i)]c(e_j)c(e_k)+\frac{1}{4}\sum_{u\neq s\neq t=1}^{n} T_{vst} a_u^v c(e_u)c(e_s)c(e_t)\nonumber\\
&-\frac{1}{2} \sum_{v,u,t=1}^{n} T_{vut} a_u^v c(e_t),\nonumber
\end{align}
where $T_{vst}:= T\big(e_v,e_s,e_t\big)$.
\end{lem}
\indent Write
 \begin{eqnarray}
{D}_x^{\alpha}&=(-i)^{|\alpha|}\partial_x^{\alpha};
~\sigma({D}_{J,T})=p_1+p_0;
~\sigma({{D}_{J,T}^{-1}})=\sum^{\infty}_{j=1}q_{-j}.\nonumber
\end{eqnarray}

\indent By the composition formula of pseudodifferential operators, we have
\begin{align}
1=\sigma({D}_{J,T}\circ {D}_{J,T}^{-1})&=\sum_{\alpha}\frac{1}{\alpha!}\partial^{\alpha}_{\xi}[\sigma({D}_{J,T})]
{D}_x^{\alpha}[\sigma({D}_{J,T}^{-1})]\nonumber\\
&=(p_1+p_0)(q_{-1}+q_{-2}+q_{-3}+\cdots)\nonumber\\
&~~~+\sum_j(\partial_{\xi_j}p_1+\partial_{\xi_j}p_0)\bigg[
({D}_{J,T})_{x_j}q_{-1}+({D}_{J,T})_{x_j}q_{-2}+({D}_{J,T})_{x_j}q_{-3}+\cdots\bigg]\nonumber\\
&=p_1q_{-1}+\bigg[p_1q_{-2}+p_0q_{-1}+\sum_j\partial_{\xi_j}p_1({D}_{J,T})_{x_j}q_{-1}\bigg]+\cdots,\nonumber
\end{align}
so
\begin{equation}
q_{-1}=p_1^{-1};~q_{-2}=-p_1^{-1}[p_0p_1^{-1}+\sum_j\partial_{\xi_j}p_1 ({D}_{J,T})_{x_j}(p_1^{-1})].\nonumber
\end{equation}
\begin{lem} The following identities hold:
\begin{align}
\sigma_{-1}({D}_{J,T}^{-1})&=\frac{ic[J(\xi)]}{|\xi|^2};\nonumber\\
\sigma_{-2}({D}_{J,T}^{-1})&=\frac{c[J(\xi)]\sigma_{0}({D}_{J,T})c[J(\xi)]}{|\xi|^4}+\frac{c[J(\xi)]}{|\xi|^6}\sum_ {j=1}^{n} c[J(dx_j)]
\Big[\partial_{x_j}(c[J(\xi)])|\xi|^2-c[J(\xi)]\partial_{x_j}(|\xi|^2)\Big].\nonumber
\end{align}
\end{lem}
\indent When $n=4$, then ${\rm tr}[{\rm \texttt{id}}]=4,$ since the sum is taken over $
r+l-k-j-|\alpha|-1=-4,~~r\leq -1,l\leq-1,$ then we have the following five cases:\\

\noindent  {\bf case a)~I)}~$r=-1,~l=-1,~k=j=0,~|\alpha|=1$\\

\noindent By applying the formula shown in \eqref{phit}, we can calculate
\begin{equation}
\Psi_1=-\int_{|\xi'|=1}\int^{+\infty}_{-\infty}\sum_{|\alpha|=1}
 {\rm tr}[\partial^\alpha_{\xi'}\pi^+_{\xi_n}\sigma_{-1}({D}_{J,T}^{-1})\times
 \partial^\alpha_{x'}\partial_{\xi_n}\sigma_{-1}({D}_{J,T}^{-1})](x_0)d\xi_n\sigma(\xi')dx'.\nonumber
\end{equation}
By Lemma 3.3, for $i<n,$ ,$J\left(d x_{p}\right)=\sum_{h=1}^{n} a_{h}^{p} d x_{h}$ then

An easy computation shows that
\begin{align}
\partial_{\xi_{n}}\partial_{x_i}\left(\frac{ic[J(\xi)]}{|\xi|^2}\right)(x_0)|_{|\xi'|=1}
=&-\frac{2i\xi_{n}}{(1+\xi_{n}^2)^2}\sum^{n}_{h=1}\sum^{n-1}_{p=1}\xi_{p}\partial_{x_i}(a^{p}_{h})c(dx_{h})\nonumber\\
&+\frac{i(1-\xi_{n}^2)}{(1+\xi_{n}^2)^2}\sum^{n}_{h=1}\partial_{x_i}(a^{n}_{h})c(dx_{h}).\nonumber
\end{align}
Likewise,
\begin{align}
\pi^+_{\xi_n}\partial_{\xi_i}\left(\frac{ic[J(\xi)]}{|\xi|^2}\right)(x_0)|_{|\xi'|=1}
&=\frac{1}{2(\xi_n-i)}c[J(dx_i)]-\frac{1}{2(\xi_n-i)^2}\xi_{i}c[J(dx_n)]\nonumber\\
&+\frac{2i-\xi_n}{2(\xi_n-i)^2}\sum^{n-1}_{q=1}\xi_{i}\xi_{q}c[J(dx_q)].\nonumber
\end{align}
On account of the above result,

We note that $i<n,$  $\int_{|\xi'|=1}{\{\xi_{i_1}\cdot\cdot\cdot\xi_{i_{2d+1}}}\}\sigma(\xi')=0,$ then we have
\begin{align}
\Psi_1
&=-\int_{|\xi'|=1}\int^{+\infty}_{-\infty}\sum_{|\alpha|=1}{\rm tr}[\partial^\alpha_{\xi'}\pi^+_{\xi_n}\sigma_{-1}({{D}_{J,T}}^{-1})\times \partial^\alpha_{x'}\partial_{\xi_n}\sigma_{-1}({{D}_{J,T}}^{-1})](x_0)d\xi_n\sigma(\xi')dx'\\
&=\int_{|\xi'|=1}\int^{+\infty}_{-\infty}\frac{i(1-\xi^2_n)}{2(\xi_n-i)^3(\xi_n+i)^2}\sum_{\beta=1}^{n}\sum_{i=1}^{n-1}a_{\beta}^{i}\partial_{x_i}(a_{\beta}^{n}){\rm tr}[\texttt{id}]d\xi_n\sigma(\xi')dx'\nonumber\\
&+\int_{|\xi'|=1}\int^{+\infty}_{-\infty}\frac{i(2i-\xi_n)(1-\xi_n^2)}{2(\xi_n-i)^4(\xi_n+i)^2}\sum_{\beta=1}^{n}\sum_{i,q=1}^{n-1}\xi_{i}\xi_{q}a_{\beta}^{q}\partial_{x_i}(a_{\beta}^{n}){\rm tr}[\texttt{id}]d\xi_n\sigma(\xi')dx'\nonumber\\
&+\int_{|\xi'|=1}\int^{+\infty}_{-\infty}\frac{i\xi_n}{(\xi_n-i)^4(\xi_n+i)^2}\sum_{\beta=1}^{n}\sum_{i,p=1}^{n-1}\xi_{i}\xi_{p}a_{\beta}^{n}\partial_{x_i}(a_{\beta}^{p}){\rm tr}[\texttt{id}]d\xi_n\sigma(\xi')dx'\nonumber\\
&=-\frac{\pi}{6}\sum_{\beta=1}^{n}\sum_{i=1}^{n-1}a_{\beta}^{i}\partial_{x_i}(a_{\beta}^{n}){\rm tr}[\texttt{id}] \Omega_3 dx' -\frac{\pi}{6}\sum_{\beta=1}^{n}\sum_{i=1}^{n-1}a_{\beta}^{n}\partial_{x_i}(a_{\beta}^{i}){\rm tr}[\texttt{id}]dx'.\nonumber
\end{align}

From \cite{Ka,LW1}, we have $\int_{|\xi'|=1}\xi_i\xi_j\sigma(\xi')=\frac{\Omega_3}{3}\delta_i^j,$ and 
$\sum_{\beta=1}^{n}\sum_{i=1}^{n-1}a_{\beta}^{i}\partial_{x_i}(a_{\beta}^{n})=\sum_{\beta,i=1}^{n}a_{\beta}^{i}\partial_{x_i}(a_{\beta}^{n})$, then
\begin{align}
\Psi_1=0.\nonumber
\end{align}

\noindent  {\bf case a)~II)}~$r=-1,~l=-1,~k=|\alpha|=0,~j=1$\\

\noindent Using \eqref{phit}, we get
\begin{equation}
\Psi_2=-\frac{1}{2}\int_{|\xi'|=1}\int^{+\infty}_{-\infty} {\rm
trace} [\partial_{x_n}\pi^+_{\xi_n}\sigma_{-1}({D}_{J,T}^{-1})\times
\partial_{\xi_n}^2\sigma_{-1}({D}_{J,T}^{-1})](x_0)d\xi_n\sigma(\xi')dx'.\nonumber
\end{equation}
It is easy to check that
\begin{align}
\pi^+_{\xi_n}\partial_{x_n}\left(\frac{ic[J(\xi)]}{|\xi|^2}\right)(x_0)|_{|\xi'|=1}=&\frac{1}{2(\xi_n-i)}\sum^{n}_{h=1}\sum^{n-1}_{p=1}\xi_{p}\partial_{x_n}(a^{p}_{h})c(dx_{h})+\frac{i}{2(\xi_n-i)}\sum^{n}_{h=1}\partial_{x_n}(a^{n}_{h})c(dx_{h})\nonumber\\
&+\frac{1}{2(\xi_n-i)}\sum^{n-1}_{p,h=1}\xi_{p}a^{p}_{h}\partial_{x_n}(c(dx_{h}))+\frac{i}{2(\xi_n-i)}\sum^{n-1}_{h=1}a^{n}_{h}\partial_{x_n}(c(dx_{h}))\nonumber\\
&+\frac{2i-\xi_n}{4(\xi_n-i)^2}h'(0)\sum^{n}_{h=1}\sum^{n-1}_{p=1}\xi_{p}a^{p}_{h}c(dx_{h})-\frac{1}{4(\xi_n-i)^2}h'(0)\sum^{n}_{h=1}a^{n}_{h}c(dx_{h}),\nonumber
\end{align}
where $\sum_{h=1}^{n-1}\partial_{x_n}(c(dx_h))=\sum_{h=1}^{n-1}\frac{1}{2}h'(0)c(dx_h).$\\
By calculation, we have
\begin{align}
\partial_{\xi_n}^2\left(\frac{ic[J(\xi)]}{|\xi|^2}\right)(x_0)|_{|\xi'|=1}
&=\frac{-2i+6i\xi_n^2}{(\xi_{n}-i)^3(\xi_{n}+i)^3}\sum^{n}_{\beta=1}\sum^{n-1}_{i=1}\xi_{i}a_{\beta}^ic(dx_{\beta})+\frac{-6i\xi_{n}+2i\xi_{n}^{3}}{(\xi_{n}-i)^3(\xi_{n}+i)^3}\sum^{n}_{\beta=1}a_{\beta}^nc(dx_{\beta}).\nonumber
\end{align}

Consequently,
\begin{align}
\Psi_2
=&-\frac{1}{2}\int_{|\xi'|=1}\int^{+\infty}_{-\infty} {\rm
tr} [\partial_{x_n}\pi^+_{\xi_n}\sigma_{-1}({D}_{J,T}^{-1})\times
\partial_{\xi_n}^2\sigma_{-1}({D}_{J,T}^{-1})](x_0)d\xi_n\sigma(\xi')dx'\\
=&\frac{\pi}{48}\sum_{\beta=1}^{n}\sum_{i=1}^{n-1}a_{\beta}^{i}\partial_{x_n}(a_{\beta}^{i}){\rm tr}[\texttt{id}]\Omega_3dx'+\frac{\pi}{16}\sum_{\beta=1}^{n}a_{\beta}^{n}\partial_{x_n}(a_{\beta}^{n}){\rm tr}[\texttt{id}]\Omega_3 dx'\nonumber\\
&+\frac{\pi}{96}\sum_{i,\beta=1}^{n-1}(a_{\beta}^{i})^2{\rm tr}[\texttt{id}]\Omega_3h'(0)dx'+\frac{\pi}{32}\sum_{\beta=1}^{n-1}(a_{\beta}^{n})^2{\rm tr}[\texttt{id}]\Omega_3h'(0)dx'\nonumber\\
&-\frac{5\pi}{192}\sum_{\beta=1}^{n}\sum_{i=1}^{n-1}(a_{\beta}^{i})^2{\rm tr}[\texttt{id}]\Omega_3h'(0)dx'-\frac{3\pi}{64}\sum_{\beta=1}^{n}(a_{\beta}^{n})^2{\rm tr}[\texttt{id}]\Omega_3h'(0)dx'.\nonumber
\end{align}

\noindent  {\bf case a)~III)}~$r=-1,~l=-1,~j=|\alpha|=0,~k=1$\\

\noindent By \eqref{phit}, we calculate that
\begin{equation}
\Psi_3=-\frac{1}{2}\int_{|\xi'|=1}\int^{+\infty}_{-\infty}
{\rm tr} [\partial_{\xi_n}\pi^+_{\xi_n}\sigma_{-1}({D}_{J,T}^{-1})\times
\partial_{\xi_n}\partial_{x_n}\sigma_{-1}({D}_{J,T}^{-1})](x_0)d\xi_n\sigma(\xi')dx'.\nonumber
\end{equation}
Hence, we have
\begin{align}
\Psi_3
=&-\frac{\pi}{48}\sum_{\beta=1}^{n}\sum_{i=1}^{n-1}a_{\beta}^{i}\partial_{x_n}(a_{\beta}^{i}){\rm tr}[\texttt{id}]\Omega_3 dx'-\frac{\pi}{96}\sum_{i,\beta=1}^{n-1}(a_{\beta}^{i})^2{\rm tr}[\texttt{id}]\Omega_3h'(0)dx'\\
&+\frac{5\pi}{192}\sum_{\beta=1}^{n}\sum_{i=1}^{n-1}(a_{\beta}^{i})^2{\rm tr}[\texttt{id}]\Omega_3h'(0)dx'
-\frac{\pi}{16}\sum_{\beta=1}^{n}a_{\beta}^{n}\partial_{x_n}(a_{\beta}^{n}){\rm tr}[\texttt{id}]\Omega_3 dx'\nonumber\\
&-\frac{\pi}{32}\sum_{\beta=1}^{n-1}(a_{\beta}^{n})^2{\rm tr}[\texttt{id}]\Omega_3h'(0)dx'\nonumber
+\frac{3\pi}{64}\sum_{\beta=1}^{n}(a_{\beta}^{n})^2{\rm tr}[\texttt{id}]\Omega_3h'(0)dx'.\nonumber
\end{align}

In combination with the calculation,
\begin{equation}
\Psi_1+\Psi_2+\Psi_3=0.\nonumber
\end{equation}

\noindent  {\bf case b)}~$r=-2,~l=-1,~k=j=|\alpha|=0$\\

\noindent Similarly, we get
\begin{align}
\Psi_4&=-i\int_{|\xi'|=1}\int^{+\infty}_{-\infty}{\rm tr} [\pi^+_{\xi_n}\sigma_{-2}({D}_{J,T}^{-1})\times
\partial_{\xi_n}\sigma_{-1}({D}_{J,T}^{-1})](x_0)d\xi_n\sigma(\xi')dx'.\nonumber
\end{align}

We first compute\\
\begin{align}
\sigma_{-2}({D}_{J,T}^{-1})(x_0)
&=\frac{c[J(\xi)]\sigma_{0}({D}_{J,T})(x_0)c[J(\xi)]}{|\xi|^4}-\frac{c[J(\xi)]}{|\xi|^6}h'(0)|\xi'|^2c[J(dx_n)]c[J(\xi)]\nonumber\\
&+\frac{c[J(\xi)]}{|\xi|^4}\Big[\sum_{j,p,h=1}^{n}\xi_p\partial_{x_j}(a_{h}^{p})c[J(dx_j)]c(dx_h)+\sum_{p=1}^{n}\sum_{h=1}^{n-1}\xi_pa_{h}^{p}c[J(dx_n)]\partial_{x_n}(c(dx_h))\Big],
\nonumber
\end{align}
where
\begin{align}
\sigma_{0}({{D}_{J,T}})(x_0)
=&-\frac{1}{4}h'(0)\sum_{\mu=1}^{n}\sum_{\nu=1}^{n-1}a_{\nu}^{\mu}c(dx_{\mu})c(dx_{n})c(dx_{\nu})+\frac{1}{4}\sum_{u\neq s\neq t,v=1}^{n} T_{vst} a_u^v c(e_u)c(e_s)c(e_t) \nonumber\\
&-\frac{1}{2} \sum_{v,u,t=1}^{n} T_{vut} a_u^v c(e_t) .\nonumber
\end{align}

For the sake of convenience in writing, we denote

\begin{align}\label{B}
	&\sigma_{-2}({D}_{J,T}^{-1})(x_0)|_{|\xi'|=1}:=B_0(x_0)+B_1(x_0)+B_2(x_0),
\end{align}

\begin{align}
B_0(x_0)=&-\frac{1}{4}h'(0)\sum_{\mu=1}^{n}\sum_{\nu=1}^{n-1}\frac{c[J(\xi)]a_{\nu}^{\mu}c(dx_{\mu})c(dx_{n})c(dx_{\nu})c[J(\xi)]}{(1+\xi_n^2)^2}-\frac{c[J(\xi)]}{(1+\xi_n^2)^3}h'(0)c[J(dx_n)]c[J(\xi)]\nonumber\\
&+\frac{c[J(\xi)]}{(1+\xi_n^2)^2}\Big[\sum_{j,p,h=1}^{n}\xi_p\partial_{x_j}(a_{h}^{p})c[J(dx_j)]c(dx_h)+\sum_{p=1}^{n}\sum_{h=1}^{n-1}\xi_pa_{h}^{p}c[J(dx_n)]\partial_{x_n}(c(dx_h))\Big];\nonumber\\
B_1(x_0)=&\ \frac{1}{4}\sum_{u\neq s\neq t,v=1}^{n} \frac{c[J(\xi)]T_{vst} a_u^v c(e_u)c(e_s)c(e_t)c[J(\xi)]}{(1+\xi_n^2)^2};\nonumber\\
B_2(x_0)=&-\frac{1}{2}\sum_{v,u,t=1}^{n}  \frac{c[J(\xi)]T_{vut} a_u^v c(e_t)c[J(\xi)]}{(1+\xi_n^2)^2},\nonumber
\end{align}
means that
\begin{align}
\pi^+_{\xi_n}\sigma_{-2}({{D}_{J,T}}^{-1})(x_0)|_{|\xi'|=1}=\pi^+_{\xi_n}(B_0(x_0))+\pi^+_{\xi_n}(B_1(x_0))+\pi^+_{\xi_n}(B_2(x_0)).\nonumber
\end{align}
By computation, we have

\begin{align}
&-i\int_{|\xi'|=1}\int^{+\infty}_{-\infty}{\rm tr} [\pi^+_{\xi_n}B_0(x_0)\times
\partial_{\xi_n}\sigma_{-1}({D}_{J,T}^{-1})](x_0)d\xi_n\sigma(\xi')dx'\\
=&-\frac{\pi}{48}\sum_{l=1}^{n}\sum_{\nu,i=1}^{n-1}((a_{\nu}^{n})^2(a_{l}^{i})^2-(a_{l}^{i})^2a_{\nu}^{\nu}a_{n}^{n}){\rm tr}[\texttt{id}]\Omega_3h'(0)dx'\nonumber\\
&-\frac{\pi}{48}\sum_{l,j,\beta=1}^{n}\sum_{i=1}^{n-1}\left((a_{\beta}^{i})^2a_{l}^{j}\partial_{x_j}(a_{l}^{n})-a_{l}^{i}a_{\beta}^{j}a_{\beta}^{i}\partial_{x_j}(a_{l}^{n})+a_{l}^{i}a_{l}^{j}a_{\beta}^{i}\partial_{x_j}(a_{\beta}^{n})\right){\rm tr}[\texttt{id}]\Omega_3 dx'\nonumber\\
&-\frac{\pi}{48}\sum_{l,j,\beta=1}^{n}\sum_{i=1}^{n-1}\left(a_{\beta}^{n}a_{l}^{j}a_{\beta}^{i}\partial_{x_j}(a_{l}^{i})-a_{l}^{n}a_{\beta}^{j}a_{\beta}^{i}\partial_{x_j}(a_{l}^{i})+a_{l}^{n}a_{l}^{j}a_{\beta}^{i}\partial_{x_j}(a_{\beta}^{i})\right){\rm tr}[\texttt{id}]\Omega_3 dx'\nonumber\\
&+\frac{\pi}{24}\sum_{l,j,\beta=1}^{n}\sum_{i=1}^{n-1}\left(a_{\beta}^{i}a_{l}^{j}a_{\beta}^{n}\partial_{x_j}(a_{l}^{i})-a_{l}^{i}a_{\beta}^{j}a_{\beta}^{n}\partial_{x_j}(a_{l}^{i})+a_{l}^{i}a_{l}^{j}a_{\beta}^{n}\partial_{x_j}(a_{\beta}^{i})\right){\rm tr}[\texttt{id}]\Omega_3 dx'\nonumber\\
&-\frac{\pi}{96}\sum_{l=1}^{n}\sum_{\nu,i=1}^{n-1}(a_{\nu}^{n})^2(a_{l}^{i})^2{\rm tr}[\texttt{id}]\Omega_3h'(0) dx'-\frac{\pi}{96}\sum_{l=1}^{n}\sum_{\nu,i=1}^{n-1}(a_{\nu}^{i})^2(a_{l}^{n})^2{\rm tr}[\texttt{id}]\Omega_3h'(0) dx'\nonumber\\
&-\frac{\pi}{48}\sum_{l=1}^{n}\sum_{\nu,i=1}^{n-1}(a_{\nu}^{i})^2(a_{l}^{n})^2{\rm tr}[\texttt{id}]\Omega_3h'(0)dx'+\frac{5\pi}{192}\sum_{\beta,l=1}^{n}\sum_{i=1}^{n-1}(a_{l}^{i})^2(a_{\beta}^{n})^2{\rm tr}[\texttt{id}]\Omega_3h'(0) dx'\nonumber\\
&-\frac{\pi}{128}\sum_{\beta,l=1}^{n}(a_{\beta}^{n})^2(a_{l}^{n})^2{\rm tr}[\texttt{id}]\Omega_3h'(0) dx'+\frac{5\pi}{128}\sum_{\beta,l=1}^{n}\sum_{i=1}^{n-1}(a_{\beta}^{i})^2(a_{l}^{n})^2{\rm tr}[\texttt{id}]\Omega_3h'(0)dx'.\nonumber
\end{align}
By computation, we have
\begin{align}\label{B_1}
\pi^+_{\xi_n}(B_1(x_0))=&\pi^+_{\xi_n}\sum_{u\neq s\neq t,v=1}^{n} \frac{c[J(\xi)]T_{vst} a_u^v c(e_u)c(e_s)c(e_t)c[J(\xi)]}{4(1+\xi_n^2)^2}\\
=&\pi^+_{\xi_n}(\frac{\xi_{n}^{2}}{4(1+\xi_{n}^{2})^{2}} \sum_{l,r,v,u\neq s\neq t=1}^{n}T_{vst} a_l^n a_r^n a_u^v c(dx_l) c(dx_u) c(dx_s) c(dx_t) c(dx_r))\nonumber\\
&+\pi^+_{\xi_n}(\frac{\xi_{n}}{4(1+\xi_{n}^{2})^{2}} \sum_{l,r,v,u\neq s\neq t=1}^{n} \sum_{\alpha=1}^{n-1} T_{vst} a_l^n a_r^{\alpha} a_u^v \xi_{\alpha} c(dx_l) c(dx_u) c(dx_s) c(dx_t) c(dx_r))\nonumber\\
&+\pi^+_{\xi_n}(\frac{\xi_{n}}{4(1+\xi_{n}^{2})^{2}} \sum_{l,r,v,u\neq s\neq t=1}^{n} \sum_{q=1}^{n-1} T_{vst} a_l^q a_r^n a_u^v \xi_{q} c(dx_l) c(dx_u) c(dx_s) c(dx_t) c(dx_r))\nonumber\\
&+\pi^+_{\xi_n}(\frac{\xi_{n}}{4(1+\xi_{n}^{2})^{2}} \sum_{l,r,v,u\neq s\neq t=1}^{n} \sum_{\alpha, q=1}^{n-1} T_{vst} a_l^q a_r^{\alpha} a_u^v \xi_{\alpha}\xi_{q} c(dx_l) c(dx_u) c(dx_s) c(dx_t) c(dx_r))\nonumber\\
=&-\frac{i \xi_{n}}{16(\xi_{n}-i)^2}\sum_{l,r,v,u\neq s\neq t=1}^{n}T_{vst} a_l^n a_r^n a_u^v c(dx_l) c(dx_u) c(dx_s) c(dx_t) c(dx_r)\nonumber\\
&-\frac{i}{16(\xi_{n}-i)^2}\sum_{l,r,v,u\neq s\neq t=1}^{n} \sum_{\alpha=1}^{n-1} T_{vst} a_l^n a_r^{\alpha} a_u^v \xi_{\alpha} c(dx_l) c(dx_u) c(dx_s) c(dx_t) c(dx_r)\nonumber\\
&-\frac{i}{16(\xi_{n}-i)^2} \sum_{l,r,v,u\neq s\neq t=1}^{n}\sum_{q=1}^{n-1} T_{vst} a_l^q a_r^n a_u^v \xi_{q} c(dx_l) c(dx_u) c(dx_s) c(dx_t) c(dx_r)\nonumber\\
&-\frac{i \xi_{n} +2}{16(\xi_{n}-i)^2}\sum_{l,r,v,u\neq s\neq t=1}^{n} \sum_{\alpha, q=1}^{n-1} T_{vst} a_l^q a_r^{\alpha} a_u^v \xi_{\alpha}\xi_{q} c(dx_l) c(dx_u) c(dx_s) c(dx_t) c(dx_r)\nonumber.
\end{align}

Thus, we have
\begin{align}
	&{\rm tr} [\pi^+_{\xi_n}B_1(x_0)\times
	\partial_{\xi_n}\sigma_{-1}({D}_{J,T}^{-1})](x_0)|_{|\xi'|=1}\\
	&=\frac{-i \xi_{n}^2}{8(\xi_{n}-i)^2(1+\xi_{n}^2)^2}\sum_{\beta, l,r,v,u\neq s\neq t=1}^{n} \sum_{i=1}^{n-1} T_{vst} \xi_{i} a_l^n a_r^n a_u^v a_{\beta}^{i} {\rm tr} [c(dx_l) c(dx_u) c(dx_s) c(dx_t) c(dx_r)c(dx_{\beta})]\nonumber\\
	&+\frac{-\xi_{n}}{8(\xi_{n}-i)^2(1+\xi_{n}^2)^2 }\sum_{\beta, l,r,v,u\neq s\neq t=1}^{n} \sum_{i,\alpha=1}^{n-1} T_{vst}\xi_{i}\xi_{\alpha} a_l^n a_r^{\alpha} a_u^v  a_{\beta}^{i} {\rm tr} [c(dx_l) c(dx_u) c(dx_s) c(dx_t) c(dx_r)c(dx_{\beta})]\nonumber\\
	&+\frac{-\xi_{n}}{8(\xi_{n}-i)^2(1+\xi_{n}^2)^2 } \sum_{\beta, l,r,v,u\neq s\neq t=1}^{n}\sum_{i,q=1}^{n-1} T_{vst} \xi_{i}\xi_{q} a_l^q a_r^n a_u^v  a_{\beta}^{i} {\rm tr} [c(dx_l) c(dx_u) c(dx_s) c(dx_t) c(dx_r)c(dx_{\beta})]\nonumber\\
	&+\frac{i \xi_{n}(i \xi_{n}+2)}{8(\xi_{n}-i)^2(1+\xi_{n}^2)^2 } \sum_{\beta, l,r,v,u\neq s\neq t=1}^{n}\sum_{i,q,\alpha=1}^{n-1} T_{vst} \xi_{i}\xi_{q} \xi_{\alpha} a_l^q a_r^{\alpha} a_u^v  a_{\beta}^{i} {\rm tr} [c(dx_l) c(dx_u) c(dx_s) c(dx_t) c(dx_r)c(dx_{\beta})]\nonumber\\
	&+\frac{\xi_{n}-\xi_{n}^{3}}{16(\xi_{n}-i)^2(1+\xi_{n}^2)^2}\sum_{\beta, l,r,v,u\neq s\neq t=1}^{n}  T_{vst}  a_l^n a_r^n a_u^v a_{\beta}^{n} {\rm tr} [c(dx_l) c(dx_u) c(dx_s) c(dx_t) c(dx_r)c(dx_{\beta})]\nonumber\\
	&+\frac{1- \xi_{n}^2}{16(\xi_{n}-i)^2(1+\xi_{n}^2)^2 }\sum_{\beta, l,r,v,u\neq s\neq t=1}^{n} \sum_{\alpha=1}^{n-1} T_{vst}\xi_{\alpha} a_l^n a_r^{\alpha} a_u^v  a_{\beta}^{n} {\rm tr} [c(dx_l) c(dx_u) c(dx_s) c(dx_t) c(dx_r)c(dx_{\beta})]\nonumber\\
	&+\frac{1-\xi_{n}^2}{16(\xi_{n}-i)^2(1+\xi_{n}^2)^2 } \sum_{\beta, l,r,v,u\neq s\neq t=1}^{n}\sum_{q=1}^{n-1} T_{vst} \xi_{q} a_l^q a_r^n a_u^v  a_{\beta}^{n} {\rm tr} [c(dx_l) c(dx_u) c(dx_s) c(dx_t) c(dx_r)c(dx_{\beta})]\nonumber\\
	&-\frac{i(1-\xi_{n}^2)(i \xi_{n}+2)}{16(\xi_{n}-i)^2(1+\xi_{n}^2)^2 } \sum_{\beta, l,r,v,u\neq s\neq t=1}^{n}\sum_{q,\alpha=1}^{n-1} T_{vst} \xi_{q} \xi_{\alpha} a_l^q a_r^{\alpha} a_u^v  a_{\beta}^{n} {\rm tr} [c(dx_l) c(dx_u) c(dx_s) c(dx_t) c(dx_r)c(dx_{\beta})]\nonumber.
\end{align}
By the relation of the Clifford action and $ {\rm tr}AB = {\rm tr}BA$, we have the equality:

\begin{align}\label{tr6}
	&\sum_{\beta, l,r,v,u\neq s\neq t=1}^{n}{\rm tr} [c(dx_l) c(dx_u) c(dx_s) c(dx_t) c(dx_r)c(dx_{\beta})]\\
	&=\sum_{\beta, l,r,v,u\neq s\neq t=1}^{n} \big[ -\delta_u^l \delta_s^{\beta} \delta_t^r+ \delta_u^l \delta_s^r \delta_t^{\beta} +\delta_u^{\beta} \delta_s^l \delta_t^r -\delta_u^{\beta} \delta_s^r \delta_t^l -\delta_u^r \delta_s^l \delta_t^{\beta} +\delta_u^r \delta_s^{\beta} \delta_t^l \big]{\rm tr}[id]\nonumber.
\end{align}

By $i<n,$  $\int_{|\xi'|=1}{\{\xi_{i_1}\cdot\cdot\cdot\xi_{i_{2d+1}}}\}\sigma(\xi')=0,\int_{\left|\xi^{\prime}\right|=1} \xi_{i} \xi_{j}\sigma(\xi')=\frac{4 \pi}{3} \delta_{i}^{j}$,then we have

\begin{align}
&-i\int_{|\xi'|=1}\int^{+\infty}_{-\infty}{\rm tr} [\pi^+_{\xi_n}B_1(x_0)\times
\partial_{\xi_n}\sigma_{-1}({{D}_{J,T}}^{-1})](x_0)d\xi_n\sigma(\xi')dx'=0.
\end{align}

By means of similar calculations, we have

\begin{align}
	&{\rm tr} [\pi^+_{\xi_n}B_2(x_0)\times
	\partial_{\xi_n}\sigma_{-1}({D}_{J,T}^{-1})](x_0)|_{|\xi'|=1}\\
	&=\frac{i \xi_{n}^2}{16(\xi_{n}-i)^2(1+\xi_{n}^2)^2}\sum_{\beta, u,v,t,l,r=1}^{n} \sum_{i=1}^{n-1} T_{vut} \xi_{i} a_l^n a_r^n a_u^v a_{\beta}^{i} {\rm tr} [c(dx_l) c(dx_t) c(dx_r)c(dx_{\beta})]\nonumber\\
    &+\frac{ \xi_{n}}{16(\xi_{n}-i)^2(1+\xi_{n}^2)^2 }\sum_{\beta, u,v,t,l,r=1}^{n} \sum_{i,\alpha=1}^{n-1} T_{vut}\xi_{i}\xi_{\alpha} a_l^n a_r^{\alpha} a_u^v  a_{\beta}^{i} {\rm tr} [c(dx_l)  c(dx_t) c(dx_r)c(dx_{\beta})]\nonumber\\
    &+\frac{ \xi_{n}}{16(\xi_{n}-i)^2(1+\xi_{n}^2)^2 } \sum_{\beta, u,v,t,l,r=1}^{n}\sum_{i,q=1}^{n-1} T_{vut} \xi_{i}\xi_{q} a_l^q a_r^n a_u^v  a_{\beta}^{i} {\rm tr} [c(dx_l) c(dx_t) c(dx_r)c(dx_{\beta})]\nonumber\\
    &-\frac{ i \xi_{n}(i \xi_{n}+2)}{16(\xi_{n}-i)^2(1+\xi_{n}^2)^2 } \sum_{\beta, u,v,t,l,r=1}^{n}\sum_{i,q,\alpha=1}^{n-1} T_{vut} \xi_{i}\xi_{q} \xi_{\alpha} a_l^q a_r^{\alpha} a_u^v  a_{\beta}^{i} {\rm tr} [c(dx_l)  c(dx_t) c(dx_r)c(dx_{\beta})]\nonumber\\
    &-\frac{\xi_{n}-\xi_{n}^{3}}{8(\xi_{n}-i)^2(1+\xi_{n}^2)^2}\sum_{\beta, u,v,t,l,r=1}^{n}  T_{vut}  a_l^n a_r^n a_u^v a_{\beta}^{n} {\rm tr} [c(dx_l) c(dx_t) c(dx_r)c(dx_{\beta})]\nonumber\\
    &-\frac{1- \xi_{n}^2}{8(\xi_{n}-i)^2(1+\xi_{n}^2)^2 }\sum_{\beta, u,v,t,l,r=1}^{n} \sum_{\alpha=1}^{n-1} T_{vut}\xi_{\alpha} a_l^n a_r^{\alpha} a_u^v  a_{\beta}^{n} {\rm tr} [c(dx_l) c(dx_t) c(dx_r)c(dx_{\beta})]\nonumber\\
    &-\frac{1-\xi_{n}^2}{8(\xi_{n}-i)^2(1+\xi_{n}^2)^2 } \sum_{\beta, u,v,t,l,r=1}^{n}\sum_{q=1}^{n-1} T_{vut} \xi_{q} a_l^q a_r^n a_u^v  a_{\beta}^{n} {\rm tr} [c(dx_l) c(dx_t) c(dx_r)c(dx_{\beta})]\nonumber\\
    &+\frac{i(1-\xi_{n}^2)(i \xi_{n}+2)}{8(\xi_{n}-i)^2(1+\xi_{n}^2)^2 } \sum_{\beta, u,v,t,l,r=1}^{n}\sum_{q,\alpha=1}^{n-1} T_{vut} \xi_{q} \xi_{\alpha} a_l^q a_r^{\alpha} a_u^v  a_{\beta}^{n} {\rm tr} [c(dx_l)  c(dx_t) c(dx_r)c(dx_{\beta})]\nonumber.
\end{align}

By the relation of the Clifford action and $ {\rm tr}AB = {\rm tr}BA$, we have the equality:

\begin{align}\label{tr4}
	&\sum_{ l,t,r,\beta=1}^{n}{\rm tr} [c(dx_l) c(dx_t) c(dx_r)c(dx_{\beta})]=\sum_{l,t,r,\beta=1}^{n} \big[ \delta_l^{\beta} \delta_t^r - \delta_l^r \delta_t^{\beta} + \delta_l^t \delta_r^{\beta} \big]{\rm tr}[id].
\end{align}
By $i<n,$  $\int_{|\xi'|=1}{\{\xi_{i_1}\cdot\cdot\cdot\xi_{i_{2d+1}}}\}\sigma(\xi')=0,\int_{\left|\xi^{\prime}\right|=1} \xi_{i} \xi_{j}\sigma(\xi')=\frac{\Omega_{3}}{3} \delta_{i}^{j}$,then we have

\begin{align}
	&-i\int_{|\xi'|=1}\int^{+\infty}_{-\infty}{\rm tr} [\pi^+_{\xi_n}B_2(x_0)\times
	\partial_{\xi_n}\sigma_{-1}({D}_{J,T}^{-1})](x_0)d\xi_n\sigma(\xi')dx'\\
	=&-i\int_{|\xi'|=1}\int^{+\infty}_{-\infty} 
	\biggl\{\frac{ \xi_{n}}{16(\xi_{n}-i)^2(1+\xi_{n}^2)^2 }\sum_{ u,v,t,r=1}^{n} \sum_{i,\alpha=1}^{n-1} T_{vut}\xi_{i}\xi_{\alpha} a_t^n a_r^{\alpha} a_u^v  a_{r}^{i}\nonumber\\
	&+\frac{ \xi_{n}}{16(\xi_{n}-i)^2(1+\xi_{n}^2)^2 }\sum_{ u,v,t,l=1}^{n} \sum_{i,q=1}^{n-1} T_{vut}\xi_{i}\xi_{q} a_q^l a_t^n a_u^v  a_{l}^{i}-\frac{\xi_{n}-\xi_{n}^{3}}{8(\xi_{n}-i)^2(1+\xi_{n}^2)^2}\sum_{u,v,t,r=1}^{n}  T_{vut}  a_t^n a_r^n a_u^v a_{r}^{n}\nonumber\\
	&-\frac{i(1-\xi_{n}^2)(i \xi_{n}+2)}{8(\xi_{n}-i)^2(1+\xi_{n}^2)^2 } \sum_{ u,v,t,l=1}^{n}\sum_{q,\alpha=1}^{n-1} T_{vut} \xi_{q} \xi_{\alpha} a_l^q a_l^{\alpha} a_u^v  a_{t}^{n}
	\biggr\}{\rm tr}[id]d\xi_n\sigma(\xi')dx'\nonumber\\
	=&\frac{\pi}{96}\sum_{ u,v,t,r=1}^{n} \sum_{i=1}^{n-1} T_{vut} a_t^n (a_r^{i})^{2} a_u^v {\Omega_{3}}{\rm tr}[id]dx'  +\frac{\pi}{96}\sum_{ u,v,t,l=1}^{n} \sum_{i=1}^{n-1} T_{vut}(a_l^i)^{2} a_t^n a_u^v{\Omega_{3}}{\rm tr}[id]dx'\nonumber\\
	&+\frac{\pi}{48}\sum_{ u,v,t,l=1}^{n}\sum_{q=1}^{n-1} T_{vut}  (a_l^q)^{2}  a_u^v  a_{t}^{n} {\Omega_{3}}{\rm tr}[id]dx'. \nonumber
\end{align}

Hence, we have
\begin{align}
	\Psi_4&=-i\int_{|\xi'|=1}\int^{+\infty}_{-\infty}{\rm tr} [\pi^+_{\xi_n}\sigma_{-2}({D}_{J,T}^{-1})\times
	\partial_{\xi_n}\sigma_{-1}({D}_{J,T}^{-1})](x_0)d\xi_n\sigma(\xi')dx'\nonumber\\
	&=-\frac{\pi}{48}\sum_{l=1}^{n}\sum_{\nu,i=1}^{n-1}((a_{\nu}^{n})^2(a_{l}^{i})^2-(a_{l}^{i})^2a_{\nu}^{\nu}a_{n}^{n}){\rm tr}[\texttt{id}]\Omega_3h'(0)dx'\nonumber\\
	&-\frac{\pi}{48}\sum_{l,j,\beta=1}^{n}\sum_{i=1}^{n-1}\left((a_{\beta}^{i})^2a_{l}^{j}\partial_{x_j}(a_{l}^{n})-a_{l}^{i}a_{\beta}^{j}a_{\beta}^{i}\partial_{x_j}(a_{l}^{n})+a_{l}^{i}a_{l}^{j}a_{\beta}^{i}\partial_{x_j}(a_{\beta}^{n})\right){\rm tr}[\texttt{id}]\Omega_3 dx'\nonumber\\
	&-\frac{\pi}{48}\sum_{l,j,\beta=1}^{n}\sum_{i=1}^{n-1}\left(a_{\beta}^{n}a_{l}^{j}a_{\beta}^{i}\partial_{x_j}(a_{l}^{i})-a_{l}^{n}a_{\beta}^{j}a_{\beta}^{i}\partial_{x_j}(a_{l}^{i})+a_{l}^{n}a_{l}^{j}a_{\beta}^{i}\partial_{x_j}(a_{\beta}^{i})\right){\rm tr}[\texttt{id}]\Omega_3 dx'\nonumber\\
	&+\frac{\pi}{24}\sum_{l,j,\beta=1}^{n}\sum_{i=1}^{n-1}\left(a_{\beta}^{i}a_{l}^{j}a_{\beta}^{n}\partial_{x_j}(a_{l}^{i})-a_{l}^{i}a_{\beta}^{j}a_{\beta}^{n}\partial_{x_j}(a_{l}^{i})+a_{l}^{i}a_{l}^{j}a_{\beta}^{n}\partial_{x_j}(a_{\beta}^{i})\right){\rm tr}[\texttt{id}]\Omega_3 dx'\nonumber\\
	&-\frac{\pi}{96}\sum_{l=1}^{n}\sum_{\nu,i=1}^{n-1}(a_{\nu}^{n})^2(a_{l}^{i})^2{\rm tr}[\texttt{id}]\Omega_3h'(0) dx'-\frac{\pi}{96}\sum_{l=1}^{n}\sum_{\nu,i=1}^{n-1}(a_{\nu}^{i})^2(a_{l}^{n})^2{\rm tr}[\texttt{id}]\Omega_3h'(0) dx'\nonumber\\
	&-\frac{\pi}{48}\sum_{l=1}^{n}\sum_{\nu,i=1}^{n-1}(a_{\nu}^{i})^2(a_{l}^{n})^2{\rm tr}[\texttt{id}]\Omega_3h'(0)dx'+\frac{5\pi}{192}\sum_{\beta,l=1}^{n}\sum_{i=1}^{n-1}(a_{l}^{i})^2(a_{\beta}^{n})^2{\rm tr}[\texttt{id}]\Omega_3h'(0) dx'\nonumber\\
	&-\frac{\pi}{128}\sum_{\beta,l=1}^{n}(a_{\beta}^{n})^2(a_{l}^{n})^2{\rm tr}[\texttt{id}]\Omega_3h'(0) dx'+\frac{5\pi}{128}\sum_{\beta,l=1}^{n}\sum_{i=1}^{n-1}(a_{\beta}^{i})^2(a_{l}^{n})^2{\rm tr}[\texttt{id}]\Omega_3h'(0)dx'\nonumber\\
	&+\frac{\pi}{96}\sum_{ u,v,t,r=1}^{n} \sum_{i=1}^{n-1} T_{vut} a_t^n (a_r^{i})^{2} a_u^v {\Omega_{3}}{\rm tr}[id]dx'  +\frac{\pi}{96}\sum_{ u,v,t,l=1}^{n} \sum_{i=1}^{n-1} T_{vut}(a_l^i)^{2} a_t^n a_u^v{\Omega_{3}}{\rm tr}[id]dx'\nonumber\\
	&+\frac{\pi}{48} \sum_{ u,v,t,l=1}^{n}\sum_{q=1}^{n-1} T_{vut}  (a_l^q)^{2}  a_u^v  a_{t}^{n} {\Omega_{3}}{\rm tr}[id]dx'.\nonumber
\end{align}

\noindent {\bf  case c)}~$r=-1,~l=-2,~k=j=|\alpha|=0$\\

\noindent We calculate
\begin{align}
\Psi_5=-i\int_{|\xi'|=1}\int^{+\infty}_{-\infty}{\rm tr} [\pi^+_{\xi_n}\sigma_{-1}({D}_{J,T}^{-1})\times
\partial_{\xi_n}\sigma_{-2}({D}_{J,T}^{-1})](x_0)d\xi_n\sigma(\xi')dx'.\nonumber
\end{align}
It is evident that
\begin{align}
\pi^+_{\xi_n}\left(\frac{ic[J(\xi)]}{|\xi|^2}\right)(x_0)|_{|\xi'|=1}
=\frac{1}{2(\xi_{n}-i)}\sum^{n}_{\beta=1}\sum^{n-1}_{i=1}\xi_{i}a_{\beta}^{i}c(dx_{\beta})+\frac{i}{2(\xi_{n}-i)}\sum^{n}_{\beta=1}a_{\beta}^{n}c(dx_{\beta}).
\end{align}
Same classification as in case b \eqref{B},
we have
\begin{align}
&-i\int_{|\xi'|=1}\int^{+\infty}_{-\infty}{\rm tr} [\pi^+_{\xi_n}\sigma_{-1}({D}_{J,T}^{-1})\times
\partial_{\xi_n}B_0](x_0)d\xi_n\sigma(\xi')dx'\\
=&\ \frac{\pi}{48}\sum_{l=1}^{n}\sum_{\nu,i=1}^{n-1}((a_{\nu}^{n})^2(a_{l}^{i})^2-(a_{l}^{i})^2a_{\nu}^{\nu}a_{n}^{n}){\rm tr}[\texttt{id}]\Omega_3h'(0)dx'\nonumber\\
&+\frac{\pi}{48}\sum_{l,j,\beta=1}^{n}\sum_{i=1}^{n-1}\left((a_{\beta}^{i})^2a_{l}^{j}\partial_{x_j}(a_{l}^{n})-a_{l}^{i}a_{\beta}^{j}a_{\beta}^{i}\partial_{x_j}(a_{l}^{n})+a_{l}^{i}a_{l}^{j}a_{\beta}^{i}\partial_{x_j}(a_{\beta}^{n})\right){\rm tr}[\texttt{id}]\Omega_3 dx'\nonumber\\
&+\frac{\pi}{48}\sum_{l,j,\beta=1}^{n}\sum_{i=1}^{n-1}\left(a_{\beta}^{n}a_{l}^{j}a_{\beta}^{i}\partial_{x_j}(a_{l}^{i})-a_{l}^{n}a_{\beta}^{j}a_{\beta}^{i}\partial_{x_j}(a_{l}^{i})+a_{l}^{n}a_{l}^{j}a_{\beta}^{i}\partial_{x_j}(a_{\beta}^{i})\right){\rm tr}[\texttt{id}]\Omega_3 dx'\nonumber\\
&-\frac{\pi}{24}\sum_{l,j,\beta=1}^{n}\sum_{i=1}^{n-1}\left(a_{\beta}^{i}a_{l}^{j}a_{\beta}^{n}\partial_{x_j}(a_{l}^{i})-a_{l}^{i}a_{\beta}^{j}a_{\beta}^{n}\partial_{x_j}(a_{l}^{i})+a_{l}^{i}a_{l}^{j}a_{\beta}^{n}\partial_{x_j}(a_{\beta}^{i})\right){\rm tr}[\texttt{id}]\Omega_3 dx'\nonumber\\
&+\frac{\pi}{96}\sum_{l=1}^{n}\sum_{\nu,i=1}^{n-1}(a_{\nu}^{n})^2(a_{l}^{i})^2{\rm tr}[\texttt{id}]\Omega_3h'(0) dx'+\frac{\pi}{96}\sum_{l=1}^{n}\sum_{\nu,i=1}^{n-1}(a_{\nu}^{i})^2(a_{l}^{n})^2{\rm tr}[\texttt{id}]\Omega_3h'(0) dx'\nonumber\\
&+\frac{\pi}{48}\sum_{l=1}^{n}\sum_{\nu,i=1}^{n-1}(a_{\nu}^{i})^2(a_{l}^{n})^2{\rm tr}[\texttt{id}]\Omega_3h'(0)dx'-\frac{5\pi}{192}\sum_{\beta,l=1}^{n}\sum_{i=1}^{n-1}(a_{l}^{i})^2(a_{\beta}^{n})^2{\rm tr}[\texttt{id}]\Omega_3h'(0) dx'\nonumber\\
&+\frac{\pi}{128}\sum_{\beta,l=1}^{n}(a_{\beta}^{n})^2(a_{l}^{n})^2{\rm tr}[\texttt{id}]\Omega_3h'(0) dx'-\frac{5\pi}{128}\sum_{\beta,l=1}^{n}\sum_{i=1}^{n-1}(a_{\beta}^{i})^2(a_{l}^{n})^2{\rm tr}[\texttt{id}]\Omega_3h'(0)dx'.\nonumber
\end{align}
By computation, we have
\begin{align}
\partial_{\xi_n}(B_1(x_0))=&\partial_{\xi_n}\sum_{u\neq s\neq t,v=1}^{n} \frac{c[J(\xi)]T_{vst} a_u^v c(e_u)c(e_s)c(e_t)c[J(\xi)]}{4(1+\xi_n^2)^2}\\
=&\partial_{\xi_n}(\frac{\xi_{n}^{2}}{4(1+\xi_{n}^{2})^{2}} \sum_{l,r,v,u\neq s\neq t=1}^{n}T_{vst} a_l^n a_r^n a_u^v c(dx_l) c(dx_u) c(dx_s) c(dx_t) c(dx_r))\nonumber\\
&+\partial_{\xi_n}(\frac{\xi_{n}}{4(1+\xi_{n}^{2})^{2}} \sum_{l,r,v,u\neq s\neq t=1}^{n} \sum_{\alpha=1}^{n-1} T_{vst} a_l^n a_r^{\alpha} a_u^v \xi_{\alpha} c(dx_l) c(dx_u) c(dx_s) c(dx_t) c(dx_r))\nonumber\\
&+\partial_{\xi_n}(\frac{\xi_{n}}{4(1+\xi_{n}^{2})^{2}} \sum_{l,r,v,u\neq s\neq t=1}^{n} \sum_{q=1}^{n-1} T_{vst} a_l^q a_r^n a_u^v \xi_{q} c(dx_l) c(dx_u) c(dx_s) c(dx_t) c(dx_r))\nonumber\\
&+\partial_{\xi_n}(\frac{\xi_{n}}{4(1+\xi_{n}^{2})^{2}} \sum_{l,r,v,u\neq s\neq t=1}^{n} \sum_{\alpha, q=1}^{n-1} T_{vst} a_l^q a_r^{\alpha} a_u^v \xi_{\alpha}\xi_{q} c(dx_l) c(dx_u) c(dx_s) c(dx_t) c(dx_r))\nonumber\\
=&\frac{\xi_{n}}{2(1+\xi_{n}^2)^2}\sum_{l,r,v,u\neq s\neq t=1}^{n}T_{vst} a_l^n a_r^n a_u^v c(dx_l) c(dx_u) c(dx_s) c(dx_t) c(dx_r)\nonumber\\
&+\frac{1-\xi_{n}^2}{4(1+\xi_{n}^2)^2}\sum_{l,r,v,u\neq s\neq t=1}^{n} \sum_{\alpha=1}^{n-1} T_{vst} a_l^n a_r^{\alpha} a_u^v \xi_{\alpha} c(dx_l) c(dx_u) c(dx_s) c(dx_t) c(dx_r)\nonumber\\
&+\frac{1-\xi_{n}^2}{4(1+\xi_{n}^2)^2} \sum_{l,r,v,u\neq s\neq t=1}^{n}\sum_{q=1}^{n-1} T_{vst} a_l^q a_r^n a_u^v \xi_{q} c(dx_l) c(dx_u) c(dx_s) c(dx_t) c(dx_r)\nonumber\\
&-\frac{\xi_{n}}{2(1+\xi_{n}^2)^2}\sum_{l,r,v,u\neq s\neq t=1}^{n} \sum_{\alpha, q=1}^{n-1} T_{vst} a_l^q a_r^{\alpha} a_u^v \xi_{\alpha}\xi_{q} c(dx_l) c(dx_u) c(dx_s) c(dx_t) c(dx_r)\nonumber.
\end{align}
Thus, we have
\begin{align}
	&{\rm tr} [\pi^+_{\xi_n}\sigma_{-1}({D}_{J,T}^{-1})\times
	\partial_{\xi_n}B_1](x_0)|_{|\xi|'=1}\\
	=&\frac{\xi_{n}}{4(\xi_{n}-i)(1+\xi_{n}^2)^2}\sum_{\beta, l,r,v,u\neq s\neq t=1}^{n} \sum_{i=1}^{n-1} T_{vst} \xi_{i} a_l^n a_r^n a_u^v a_{\beta}^{i} {\rm tr} [c(dx_l) c(dx_u) c(dx_s) c(dx_t) c(dx_r)c(dx_{\beta})]\nonumber\\
	&+\frac{1-\xi_{n}^2}{8(\xi_{n}-i)(1+\xi_{n}^2)^2 }\sum_{\beta, l,r,v,u\neq s\neq t=1}^{n} \sum_{i,\alpha=1}^{n-1} T_{vst}\xi_{i}\xi_{\alpha} a_l^n a_r^{\alpha} a_u^v  a_{\beta}^{i} {\rm tr} [c(dx_l) c(dx_u) c(dx_s) c(dx_t) c(dx_r)c(dx_{\beta})]\nonumber\\
	&+\frac{1-\xi_{n}^2}{8(\xi_{n}-i)(1+\xi_{n}^2)^2} \sum_{\beta, l,r,v,u\neq s\neq t=1}^{n}\sum_{i,q=1}^{n-1} T_{vst} \xi_{i}\xi_{q} a_l^q a_r^n a_u^v  a_{\beta}^{i} {\rm tr} [c(dx_l) c(dx_u) c(dx_s) c(dx_t) c(dx_r)c(dx_{\beta})]\nonumber\\
	&-\frac{\xi_{n}}{4(\xi_{n}-i)(1+\xi_{n}^2)^2} \sum_{\beta, l,r,v,u\neq s\neq t=1}^{n}\sum_{i,q,\alpha=1}^{n-1} T_{vst} \xi_{i}\xi_{q} \xi_{\alpha} a_l^q a_r^{\alpha} a_u^v  a_{\beta}^{i} {\rm tr} [c(dx_l) c(dx_u) c(dx_s) c(dx_t) c(dx_r)c(dx_{\beta})]\nonumber\\
	&+\frac{i\xi_{n}}{4(\xi_{n}-i)(1+\xi_{n}^2)^2}\sum_{\beta, l,r,v,u\neq s\neq t=1}^{n}  T_{vst}  a_l^n a_r^n a_u^v a_{\beta}^{n} {\rm tr} [c(dx_l) c(dx_u) c(dx_s) c(dx_t) c(dx_r)c(dx_{\beta})]\nonumber\\
	&+\frac{i- i\xi_{n}^2}{8(\xi_{n}-i)(1+\xi_{n}^2)^2 }\sum_{\beta, l,r,v,u\neq s\neq t=1}^{n} \sum_{\alpha=1}^{n-1} T_{vst}\xi_{\alpha} a_l^n a_r^{\alpha} a_u^v  a_{\beta}^{n} {\rm tr} [c(dx_l) c(dx_u) c(dx_s) c(dx_t) c(dx_r)c(dx_{\beta})]\nonumber\\
	&+\frac{i-i\xi_{n}^2}{8(\xi_{n}-i)(1+\xi_{n}^2)^2 } \sum_{\beta, l,r,v,u\neq s\neq t=1}^{n}\sum_{q=1}^{n-1} T_{vst} \xi_{q} a_l^q a_r^n a_u^v  a_{\beta}^{n} {\rm tr} [c(dx_l) c(dx_u) c(dx_s) c(dx_t) c(dx_r)c(dx_{\beta})]\nonumber\\
	&-\frac{i \xi_{n}}{8(\xi_{n}-i)(1+\xi_{n}^2)^2 } \sum_{\beta, l,r,v,u\neq s\neq t=1}^{n}\sum_{q,\alpha=1}^{n-1} T_{vst} \xi_{q} \xi_{\alpha} a_l^q a_r^{\alpha} a_u^v  a_{\beta}^{n} {\rm tr} [c(dx_l) c(dx_u) c(dx_s) c(dx_t) c(dx_r)c(dx_{\beta})]\nonumber.
\end{align}
According to \eqref{tr6}, we have
\begin{align}
-i\int_{|\xi'|=1}\int^{+\infty}_{-\infty}{\rm tr} [\pi^+_{\xi_n}\sigma_{-1}({D}_{J,T}^{-1})\times
\partial_{\xi_n}B_1](x_0)d\xi_n\sigma(\xi')dx'=0
\end{align}

By means of similar calculations, we have
\begin{align}
	\partial_{\xi_n}(B_2(x_0))=&\partial_{\xi_n}\sum_{u=s,v,t=1}^{n} -\frac{c[J(\xi)]T_{vut} a_u^v c(e_t)c[J(\xi)]}{2(1+\xi_n^2)^2}\\
	=&-\frac{\xi_{n}}{(1+\xi_{n}^2)^2}\sum_{l,r,v,u, t=1}^{n}T_{vut} a_l^n a_r^n a_u^v c(dx_l)  c(dx_t) c(dx_r)\nonumber\\
	&-\frac{1-\xi_{n}^2}{2(1+\xi_{n}^2)^2}\sum_{l,r,v,u, t=1}^{n} \sum_{\alpha=1}^{n-1} T_{vut} a_l^n a_r^{\alpha} a_u^v \xi_{\alpha} c(dx_l) c(dx_t) c(dx_r)\nonumber\\
	&-\frac{1-\xi_{n}^2}{2(1+\xi_{n}^2)^2} \sum_{l,r,v,u, t=1}^{n}\sum_{q=1}^{n-1} T_{vut} a_l^q a_r^n a_u^v \xi_{q} c(dx_l) c(dx_t) c(dx_r)\nonumber\\
	&+\frac{\xi_{n}}{(1+\xi_{n}^2)^2}\sum_{l,r,v,u, t=1}^{n} \sum_{\alpha, q=1}^{n-1} T_{vut} a_l^q a_r^{\alpha} a_u^v \xi_{\alpha}\xi_{q} c(dx_l)  c(dx_t) c(dx_r)\nonumber.
\end{align}
Then
\begin{align}
	&{\rm tr} [\pi^+_{\xi_n}B_2(x_0)\times
	\partial_{\xi_n}\sigma_{-1}({D}_{J,T}^{-1})](x_0)|_{|\xi|'=1}\\
	=&-\frac{\xi_{n}}{2(\xi_{n}-i)(1+\xi_{n}^2)^2}\sum_{\beta, u,v,t,l,r=1}^{n} \sum_{i=1}^{n-1} T_{vut} \xi_{i} a_l^n a_r^n a_u^v a_{\beta}^{i} {\rm tr} [c(dx_l) c(dx_t) c(dx_r)c(dx_{\beta})]\nonumber\\
	&-\frac{1-\xi_{n}^2}{4(\xi_{n}-i)(1+\xi_{n}^2)^2}\sum_{\beta, u,v,t,l,r=1}^{n} \sum_{i,\alpha=1}^{n-1} T_{vut}\xi_{i}\xi_{\alpha} a_l^n a_r^{\alpha} a_u^v  a_{\beta}^{i} {\rm tr} [c(dx_l)  c(dx_t) c(dx_r)c(dx_{\beta})]\nonumber\\
	&-\frac{1-\xi_{n}^2}{4(\xi_{n}-i)(1+\xi_{n}^2)^2} \sum_{\beta, u,v,t,l,r=1}^{n}\sum_{i,q=1}^{n-1} T_{vut} \xi_{i}\xi_{q} a_l^q a_r^n a_u^v  a_{\beta}^{i} {\rm tr} [c(dx_l) c(dx_t) c(dx_r)c(dx_{\beta})]\nonumber\\
	&+\frac{\xi_{n}}{2(\xi_{n}-i)(1+\xi_{n}^2)^2} \sum_{\beta, u,v,t,l,r=1}^{n}\sum_{i,q,\alpha=1}^{n-1} T_{vut} \xi_{i}\xi_{q} \xi_{\alpha} a_l^q a_r^{\alpha} a_u^v  a_{\beta}^{i} {\rm tr} [c(dx_l)  c(dx_t) c(dx_r)c(dx_{\beta})]\nonumber\\
	&-\frac{i\xi_{n}}{2(\xi_{n}-i)(1+\xi_{n}^2)^2}\sum_{\beta, u,v,t,l,r=1}^{n}  T_{vut}  a_l^n a_r^n a_u^v a_{\beta}^{n} {\rm tr} [c(dx_l) c(dx_t) c(dx_r)c(dx_{\beta})]\nonumber\\
	&-\frac{i-i\xi_{n}}{4(\xi_{n}-i)(1+\xi_{n}^2)^2}\sum_{\beta, u,v,t,l,r=1}^{n} \sum_{\alpha=1}^{n-1} T_{vut}\xi_{\alpha} a_l^n a_r^{\alpha} a_u^v  a_{\beta}^{n} {\rm tr} [c(dx_l) c(dx_t) c(dx_r)c(dx_{\beta})]\nonumber\\
	&-\frac{i-i\xi_{n}}{4(\xi_{n}-i)(1+\xi_{n}^2)^2} \sum_{\beta, u,v,t,l,r=1}^{n}\sum_{q=1}^{n-1} T_{vut} \xi_{q} a_l^q a_r^n a_u^v  a_{\beta}^{n} {\rm tr} [c(dx_l) c(dx_t) c(dx_r)c(dx_{\beta})]\nonumber\\
	&+\frac{\xi_{n}}{2(\xi_{n}-i)(1+\xi_{n}^2)^2} \sum_{\beta, u,v,t,l,r=1}^{n}\sum_{q,\alpha=1}^{n-1} T_{vut} \xi_{q} \xi_{\alpha} a_l^q a_r^{\alpha} a_u^v  a_{\beta}^{n} {\rm tr} [c(dx_l)  c(dx_t) c(dx_r)c(dx_{\beta})]\nonumber.
\end{align}
According to \eqref{tr4}, we have
\begin{align}
	&-i\int_{|\xi'|=1}\int^{+\infty}_{-\infty}{\rm tr} [\pi^+_{\xi_n}\sigma_{-1}({D}_{J,T}^{-1})\times
	\partial_{\xi_n}B_2](x_0)d\xi_n\sigma(\xi')dx'\\
	=&-i\int_{|\xi'|=1}\int^{+\infty}_{-\infty} \biggl\{-\frac{1-\xi_{n}^2}{4(\xi_{n}-i)^3(\xi_{n}+i)^2}\sum_{ u,v,t,r=1}^{n} \sum_{i,\alpha=1}^{n-1} T_{vut}\xi_{i}\xi_{\alpha} a_t^n a_r^{\alpha} a_u^v  a_r^{i}\nonumber\\
	&-\frac{1-\xi_{n}^2}{4(\xi_{n}-i)^3(\xi_{n}+i)^2} \sum_{ u,v,l,r=1}^{n}\sum_{i,q=1}^{n-1} T_{vut} \xi_{i}\xi_{q} a_l^q a_t^n a_u^v a_l^i \nonumber\\
	&-\frac{i\xi_{n}}{2(\xi_{n}-i)^3(\xi_{n}+i)^2}\sum_{ u,v,t,r=1}^{n}  T_{vut}  a_t^n (a_r^n)^2 a_u^v \nonumber\\
	&-\frac{i\xi_{n}}{2(\xi_{n}-i)^3(\xi_{n}+i)^2} \sum_{ u,v,t,l=1}^{n}\sum_{q,\alpha=1}^{n-1} T_{vut} \xi_{q} \xi_{\alpha} a_l^q a_l^{\alpha} a_u^v  a_t^{n} \biggr\}{\rm tr} [id]d\xi_n\sigma(\xi')dx'\nonumber\\
	=&-\frac{\pi}{48}\sum_{ u,v,t,r=1}^{n} \sum_{i=1}^{n-1} T_{vut} a_t^n a_u^v  (a_r^{i})^2{\rm tr} [id]\Omega_{3}dx'-\frac{\pi}{48}\sum_{ u,v,l,r=1}^{n}\sum_{i=1}^{n-1} T_{vut} (a_l^i)^2 a_t^n a_u^v {\rm tr} [id]\Omega_{3}dx' \nonumber\\
	&-\frac{\pi}{16}\sum_{ u,v,t,r=1}^{n}  T_{vut}  a_t^n (a_r^n)^2 a_u^v {\rm tr} [id]\Omega_{3}dx'-\frac{\pi}{48}\sum_{ u,v,t,l=1}^{n}\sum_{i=1}^{n-1} T_{vut} (a_l^i)^2  a_u^v  a_t^{n} {\rm tr} [id]\Omega_{3}dx'\nonumber.
\end{align}

Consequently
\begin{align}
	\Psi_5&=-i\int_{|\xi'|=1}\int^{+\infty}_{-\infty}{\rm tr} [\pi^+_{\xi_n}\sigma_{-1}({D}_{J,T}^{-1})\times
	\partial_{\xi_n}\sigma_{-2}({D}_{J,T}^{-1})](x_0)d\xi_n\sigma(\xi')dx'\\
	&=\ \frac{\pi}{48}\sum_{l=1}^{n}\sum_{\nu,i=1}^{n-1}((a_{\nu}^{n})^2(a_{l}^{i})^2-(a_{l}^{i})^2a_{\nu}^{\nu}a_{n}^{n}){\rm tr}[\texttt{id}]\Omega_3h'(0)dx'\nonumber\\
	&+\frac{\pi}{48}\sum_{l,j,\beta=1}^{n}\sum_{i=1}^{n-1}\left((a_{\beta}^{i})^2a_{l}^{j}\partial_{x_j}(a_{l}^{n})-a_{l}^{i}a_{\beta}^{j}a_{\beta}^{i}\partial_{x_j}(a_{l}^{n})+a_{l}^{i}a_{l}^{j}a_{\beta}^{i}\partial_{x_j}(a_{\beta}^{n})\right){\rm tr}[\texttt{id}]\Omega_3 dx'\nonumber\\
	&+\frac{\pi}{48}\sum_{l,j,\beta=1}^{n}\sum_{i=1}^{n-1}\left(a_{\beta}^{n}a_{l}^{j}a_{\beta}^{i}\partial_{x_j}(a_{l}^{i})-a_{l}^{n}a_{\beta}^{j}a_{\beta}^{i}\partial_{x_j}(a_{l}^{i})+a_{l}^{n}a_{l}^{j}a_{\beta}^{i}\partial_{x_j}(a_{\beta}^{i})\right){\rm tr}[\texttt{id}]\Omega_3 dx'\nonumber\\
	&-\frac{\pi}{24}\sum_{l,j,\beta=1}^{n}\sum_{i=1}^{n-1}\left(a_{\beta}^{i}a_{l}^{j}a_{\beta}^{n}\partial_{x_j}(a_{l}^{i})-a_{l}^{i}a_{\beta}^{j}a_{\beta}^{n}\partial_{x_j}(a_{l}^{i})+a_{l}^{i}a_{l}^{j}a_{\beta}^{n}\partial_{x_j}(a_{\beta}^{i})\right){\rm tr}[\texttt{id}]\Omega_3 dx'\nonumber\\
	&+\frac{\pi}{96}\sum_{l=1}^{n}\sum_{\nu,i=1}^{n-1}(a_{\nu}^{n})^2(a_{l}^{i})^2{\rm tr}[\texttt{id}]\Omega_3h'(0) dx'+\frac{\pi}{96}\sum_{l=1}^{n}\sum_{\nu,i=1}^{n-1}(a_{\nu}^{i})^2(a_{l}^{n})^2{\rm tr}[\texttt{id}]\Omega_3h'(0) dx'\nonumber\\
	&+\frac{\pi}{48}\sum_{l=1}^{n}\sum_{\nu,i=1}^{n-1}(a_{\nu}^{i})^2(a_{l}^{n})^2{\rm tr}[\texttt{id}]\Omega_3h'(0)dx'-\frac{5\pi}{192}\sum_{\beta,l=1}^{n}\sum_{i=1}^{n-1}(a_{l}^{i})^2(a_{\beta}^{n})^2{\rm tr}[\texttt{id}]\Omega_3h'(0) dx'\nonumber\\
	&+\frac{\pi}{128}\sum_{\beta,l=1}^{n}(a_{\beta}^{n})^2(a_{l}^{n})^2{\rm tr}[\texttt{id}]\Omega_3h'(0) dx'-\frac{5\pi}{128}\sum_{\beta,l=1}^{n}\sum_{i=1}^{n-1}(a_{\beta}^{i})^2(a_{l}^{n})^2{\rm tr}[\texttt{id}]\Omega_3h'(0)dx'\nonumber\\
	&-\frac{\pi}{48}\sum_{ u,v,t,r=1}^{n} \sum_{i=1}^{n-1} T_{vut} a_t^n a_u^v  (a_r^{i})^2{\rm tr} [id]\Omega_{3}dx'-\frac{\pi}{48}\sum_{ u,v,l,r=1}^{n}\sum_{i=1}^{n-1} T_{vut} (a_l^i)^2 a_t^n a_u^v {\rm tr} [id]\Omega_{3}dx' \nonumber\\
	&-\frac{\pi}{16}\sum_{ u,v,t,r=1}^{n}  T_{vut}  a_t^n (a_r^n)^2 a_u^v {\rm tr} [id]\Omega_{3}dx'-\frac{\pi}{48}\sum_{ u,v,t,l=1}^{n}\sum_{i=1}^{n-1} T_{vut} (a_l^i)^2  a_u^v  a_t^{n} {\rm tr} [id]\Omega_{3}dx'\nonumber.
\end{align}

In summary,
\begin{align}
\Psi=\Psi_1+\Psi_2+\Psi_3+\Psi_4+\Psi_5=0.\nonumber
\end{align}

Thus, the following theorem is obtained.
\begin{thm}
Let $M$ be a $4$-dimensional almost product Riemannian spin manifold with the boundary $\partial M$ and the metric
$g^M$ as above, ${{D}_{J,T}}$ be the $J$-twist of the Dirac operator on $\widetilde{M}$, then
\begin{align}
&\widetilde{{\rm Wres}}[\pi^+{D}_{J,T}^{-1}\circ\pi^+{D}_{J,T}^{-1}]\\
=&\int_{M}2\pi^{2}\biggl\{ \sum_{i,j=1}^{n}R(J(e_{i}), J(e_{j}), e_{j}, e_{i})
-2\sum_{\nu,j=1}^{n}g^{M}(\nabla_{e_{j}}^{L}(J)e_{\nu}, (\nabla^{L}_{e_{\nu}}J)e_{j})\nonumber \\
&-2\sum_{\nu,j=1}^{n}g^{M}(J(e_{\nu}), (\nabla^{L}_{e_{j}}(\nabla^{L}_{e_{\nu}}(J)))e_{j}-(\nabla^{L}_{\nabla^{L}_{e_{j}}e_{\nu}}(J))e_{j})\nonumber\\
&-\sum_{\alpha,\nu,j=1}^{n}g^{M}(J(e_{\alpha}), (\nabla^{L}_{e_{\nu}}J)e_{j})g^{M}((\nabla^{L}_{e_{\alpha}}J)e_{j}, J(e_{\nu}))\nonumber\\
	&-\sum_{\alpha,\nu,j=1}^{n}g^{M}(J(e_{\alpha}), (\nabla^{L}_{e_{\alpha}}J)e_{j})g^{M}(J(e_{\nu}), (\nabla^{L}_{e_{\nu}}J)e_{j})+\sum_{\nu,j=1}^{n}g^{M}((\nabla^{L}_{e_{\nu}}J)e_{j}, (\nabla^{L}_{e_{\nu}}J)e_{j})-\frac{5}{3}s\nonumber\\
&-2\sum_{j, m, p} T\left(e_m, e_{p}, J\left(e_{m}\right)\right) g^{M}\left((\nabla_{e_{j}}^{L}J) e_{j} ,e_{p}\right)-2\sum_{j, m, \alpha} T\left(e_m, J\left(e_{j}\right), J\left(e_{\alpha}\right)\right) g^{M}\left((\nabla_{e_{\alpha}}^{L}J) e_{j} ,J(e_m)\right)\nonumber\\
&-2\sum_{j, m, \alpha} T\left(e_m, J\left(e_{j}\right), J\left(e_{m}\right)\right) g^{M}\left((\nabla_{e_{\alpha}}^{L}J) e_{j} ,J(e_{\alpha})\right)+2\sum_{m,p, \alpha} T\left(e_m, e_p, J\left(e_{\alpha}\right)\right) g^{M}\left((\nabla_{e_{\alpha}}^{L}J) e_{m} ,e_p\right)\nonumber\\
&-2\sum_{m,p, \alpha} T\left(e_m, e_p, J\left(e_{\alpha}\right)\right) g^{M}\left((\nabla_{e_{m}}^{L}J) e_{\alpha} ,e_p\right)- \sum_{j,p,m} T(e_j, e_p, J(e_m))T(e_m, e_p, J(e_j))\nonumber\\
&-2 \sum_{j, l, \alpha} T^2\left(J\left(e_{j}\right), e_{l}, e_{\alpha}\right)\biggr\}d{\rm Vol_M}\nonumber
\end{align}
where $s$ is the scalar curvature.
\end{thm}

\section{ A Kastler-Kalau-Walze type theorem for $6$-dimensional manifolds with boundary}

\begin{lem} The following identities hold:
	\begin{align}
		\sigma_3({D}_{J,T}^3)=&\sqrt{-1}c[J(\xi)]|\xi|^2;\nonumber\\
		\sigma_2({D}_{J,T}^3)=
		&\sum_{i, j, l=1}^{n} c\left[J\left(d x_{l}\right)\right] \partial_{l}\left(g^{i j}\right) \xi_{i} \xi_{j}+c[J(\xi)]\left(2\sigma^{k}-\Gamma^{k}\right) \xi_{k}-\sum_{\alpha=1}^{n} c[J(\xi)] c\left[J\left(e_{\alpha}\right)\right] c\left[\left(\nabla_{e_{\alpha}}^{L} J\right)\left(\xi^{*}\right)\right] \nonumber\\
		&-\frac{1}{4}|\xi|^{2} \sum_{s, t, l=1}^{n} \omega_{s, t}\left(e_{l}\right) c\left[J\left(e_{l}\right)\right] c\left(e_{s}\right) c\left(e_{t}\right)-\sum_{j, l=1}^{n}c[J(e_l)]c[J(e_j)]T_J \xi_{l}\xi_{j}\nonumber\\
		&-\sum_{j, l=1}^{n}c[J(e_l)]T_Jc[J(e_j)] \xi_{l}\xi_{j}+T_J |\xi|^2,\nonumber
	\end{align}
where $\xi^{*}=\sum_{\beta=1}^{n}\left\langle e_{\beta},\xi \right\rangle e_{\beta} .$
\end{lem}
\indent Write
\begin{eqnarray}
	D_x^{\alpha}&=(-i)^{|\alpha|}\partial_x^{\alpha};
	~\sigma({D}_{J,T}^3)=p_3+p_2+p_1+p_0;
	~\sigma({{D}_{J,T}^{-3}})=\sum^{\infty}_{j=3}q_{-j}.\nonumber
\end{eqnarray}

\indent By the composition formula of pseudodifferential operators, we have
\begin{align}
	1=\sigma({D}_{J,T}\circ {D}_{J,T}^{-1})&=\sum_{\alpha}\frac{1}{\alpha!}\partial^{\alpha}_{\xi}[\sigma({D}_{J,T})]
	{D}_x^{\alpha}[\sigma({D}_{J,T}^{-1})]\nonumber\\
	&=(p_3+p_2+p_1+p_0)(q_{-3}+q_{-4}+q_{-5}+\cdots)\nonumber\\
	&~~~+\sum_j(\partial_{\xi_j}p_3+\partial_{\xi_j}p_2+\partial_{\xi_j}p_1+\partial_{\xi_j}p_0)(
	D_{x_j}q_{-3}+D_{x_j}q_{-4}+D_{x_j}q_{-5}+\cdots)\nonumber\\
	&=p_3q_{-3}+(p_3q_{-4}+p_2q_{-3}+\sum_j\partial_{\xi_j}p_3D_{x_j}q_{-3})+\cdots,\nonumber
\end{align}
so
\begin{equation}
	q_{-3}=p_3^{-1};~q_{-4}=-p_3^{-1}[p_2p_3^{-1}+\sum_j\partial_{\xi_j}p_3D_{x_j}(p_3^{-1})].\nonumber
\end{equation}
\begin{lem} The following identities hold:
	\begin{align}
		\sigma_{-3}({D}_{J,T}^{-3})&=\frac{ic[J(\xi)]}{|\xi|^4};\nonumber\\
		\sigma_{-4}({{D}_{J,T}^{-3}})&=\frac{c[J(\xi)]\sigma_{2}({D}_{J,T}^3)c[J(\xi)]}{|\xi|^8}\nonumber\\
		&\;\;\;+\frac{c[J(\xi)]}{|\xi|^{10}}\sum_ {j=1}^{n} \Big[c[J(dx_j)|\xi|^2+2\xi_{j}c[J(\xi)]\Big]
		\Big[\partial_{x_j}(c[J(\xi)])|\xi|^2-c[J(\xi)]\partial_{x_j}(|\xi|^2)\Big].\nonumber
	\end{align}
\end{lem}
\indent When $n=6$, then ${\rm tr}[{\rm \texttt{id}}]=8,$ since the sum is taken over $
r+l-k-j-|\alpha|-1=-6,~~r\leq -1,l\leq-3,$ then we have the following five cases:\\

\noindent  {\bf case $\bf \hat{a}$)~I)}~$r=-1,~l=-3,~k=j=0,~|\alpha|=1$\\

\noindent By applying the formula shown in \eqref{phi}, we can calculate
\begin{equation}
	\hat{\Psi}_1=-\int_{|\xi'|=1}\int^{+\infty}_{-\infty}\sum_{|\alpha|=1}
	{\rm tr}[\partial^\alpha_{\xi'}\pi^+_{\xi_n}\sigma_{-1}({D}_{J,T}^{-1})\times
	\partial^\alpha_{x'}\partial_{\xi_n}\sigma_{-3}({D}_{J,T}^{-3})](x_0)d\xi_n\sigma(\xi')dx'.\nonumber
\end{equation}
By Lemma 3.3 and Lemma 4.2, for $i<n,$ $J\left(d x_{p}\right)=\sum_{h=1}^{n} a_{h}^{p} d x_{h}$. Then by a simple calculation we get
\begin{align}
	\hat{\Psi_1}=0.\nonumber
\end{align}

\noindent  {\bf case $\bf \hat{a}$)~II)}~$r=-1,~l=-3,~k=|\alpha|=0,~j=1$\\

\noindent Using \eqref{phi}, we get
\begin{equation}
	\hat{\Psi_2}=-\frac{1}{2}\int_{|\xi'|=1}\int^{+\infty}_{-\infty} {\rm
		trace} [\partial_{x_n}\pi^+_{\xi_n}\sigma_{-1}({D}_{J,T}^{-1})\times
	\partial_{\xi_n}^2\sigma_{-3}({D}_{J,T}^{-3})](x_0)d\xi_n\sigma(\xi')dx'.\nonumber
\end{equation}
A simple calculation yields
\begin{align}
	\hat{\Psi_2}
	=&\;\;\;\frac{\pi}{64}\sum_{h=1}^{n} \sum_{i=1}^{n-1} a_{h}^{i} \partial_{x_{n}}\left(a_{h}^{i}\right) \operatorname{tr}[\mathrm{id}] \Omega_{4} d x^{\prime}+\frac{3 \pi}{64}\sum_{h=1}^{n} a_{h}^{n} \partial_{x_{n}}\left(a_{h}^{n}\right) \operatorname{tr}[\mathrm{id}] \Omega_{4} d x^{\prime} \nonumber \\
	&+\frac{\pi}{128}\sum_{h, i=1}^{n-1}\left(a_{h}^{i}\right)^{2} \operatorname{tr}[\mathrm{id}] \Omega_{4} h^{\prime}(0)d x^{\prime}+\frac{3 \pi}{128}\sum_{h=1}^{n-1}\left(a_{h}^{n}\right)^{2} \operatorname{tr}[\mathrm{id}] \Omega_{4} h^{\prime}(0)d x^{\prime}\nonumber \\
	&-\frac{7 \pi}{320}\sum_{h=1}^{n} \sum_{i=1}^{n-1}\left(a_{h}^{i}\right)^{2} \operatorname{tr}[\mathrm{id}] \Omega_{4} h^{\prime}(0) d x^{\prime}-\frac{3 \pi}{64}\sum_{h=1}^{n}\left(a_{h}^{n}\right)^{2} \operatorname{tr}[\mathrm{id}] \Omega_{4} h^{\prime}(0) d x^{\prime}.\nonumber
\end{align}

\noindent  {\bf case $\bf \hat{a}$)~III)}~$r=-1,~l=-3,~j=|\alpha|=0,~k=1$\\

\noindent By \eqref{phi}, we calculate that
\begin{equation}
	\hat{\Psi_3}=-\frac{1}{2}\int_{|\xi'|=1}\int^{+\infty}_{-\infty}
	{\rm tr} [\partial_{\xi_n}\pi^+_{\xi_n}\sigma_{-1}({D}_{J,T}^{-1})\times
	\partial_{\xi_n}\partial_{x_n}\sigma_{-3}({D}_{J,T}^{-3})](x_0)d\xi_n\sigma(\xi')dx'.\nonumber\\
\end{equation}
Hence, we have
\begin{align}
	\hat{\Psi_3}
=&-\frac{\pi}{64}\sum_{h=1}^{n} \sum_{i=1}^{n-1} a_{h}^{i} \partial_{x_{n}}\left(a_{h}^{i}\right) \operatorname{tr}[\mathrm{id}] \Omega_{4} d x^{\prime}-\frac{\pi}{128}\sum_{i, h=1}^{n-1}\left(a_{h}^{i}\right)^{2} \operatorname{tr}[\mathrm{id}] \Omega_{4} h^{\prime}(0) d x^{\prime} \nonumber \\
&+\frac{21 \pi}{640}\sum_{h=1}^{n} \sum_{i=1}^{n-1}\left(a_{h}^{i}\right)^{2} \operatorname{tr}[\mathrm{id}] \Omega_{4} h^{\prime}(0) d x^{\prime}-\frac{3 \pi}{64}\sum_{h=1}^{n} a_{h}^{n} \partial_{x_{n}}\left(a_{h}^{n}\right) \operatorname{tr}[\mathrm{id}] \Omega_{4} d x^{\prime}\nonumber \\
&-\frac{3 \pi}{128}\sum_{h=1}^{n-1}\left(a_{h}^{n}\right)^{2} \operatorname{tr}[\mathrm{id}] \Omega_{4} h^{\prime}(0) d x^{\prime}+\frac{9 \pi}{128}\sum_{h=1}^{n}\left(a_{h}^{n}\right)^{2} \operatorname{tr}[\mathrm{id}] \Omega_{4} h^{\prime}(0)d x^{\prime} .\nonumber
\end{align}

In combination with the calculation,
\begin{align}
	\hat{\Psi_1}+\hat{\Psi_2}+\hat{\Psi_3} =\;\frac{7 \pi}{640}\sum_{l=1}^{n} \sum_{i=1}^{n-1}\left(a_{l}^{i}\right)^{2} \operatorname{tr}[\mathrm{id}] \Omega_{4} h^{\prime}(0) d x^{\prime} +\frac{3 \pi}{128}\sum_{l=1}^{n}\left(a_{l}^{n}\right)^{2} \operatorname{tr}[\mathrm{id}] \Omega_{4} h^{\prime}(0) d x^{\prime} .\nonumber
\end{align}

\noindent  { \bf case $\bf \hat{b}$)}~$r=-1,~l=-4,~k=j=|\alpha|=0$\\
\noindent By \eqref{phi}, we calculate that
\begin{align}
	\hat{\Psi_4}=&-i\int_{|\xi'|=1}\int^{+\infty}_{-\infty}
	{\rm tr} [\pi^+_{\xi_n}\sigma_{-1}({D}_{J,T}^{-1})\times
	\partial_{\xi_n}\sigma_{-3}({D}_{J,T}^{-4})](x_0)d\xi_n\sigma(\xi')dx'\nonumber\\
	=&\;i\int_{|\xi'|=1}\int^{+\infty}_{-\infty}
	{\rm tr} [\partial_{\xi_n}\pi^+_{\xi_n}\sigma_{-1}({D}_{J,T}^{-1})\times \sigma_{-3}({D}_{J,T}^{-4})](x_0)d\xi_n\sigma(\xi')dx'.\nonumber
\end{align}

From Lemma 4.1 and Lemma 4.2 we get

\begin{align}
	\sigma_{-4}({{D}_{J,T}^{-3}})=&\;\;\;\frac{c[J(\xi)]\sigma_{2}({D}_{J,T}^3)c[J(\xi)]}{|\xi|^8}\nonumber\\
    &+\frac{c[J(\xi)]}{|\xi|^{10}}\sum_ {j=1}^{n} \Big[c[J(dx_j)|\xi|^2+2\xi_{j}c[J(\xi)]\Big]
    \Big[\partial_{x_j}(c[J(\xi)])|\xi|^2-c[J(\xi)]\partial_{x_j}(|\xi|^2)\Big]\nonumber\\
    =&\;\;\;\frac{c[J(\xi)]\sigma_{2}({D}_{J}^3)c[J(\xi)]}{|\xi|^8}\nonumber\\
    &+\frac{c[J(\xi)]}{|\xi|^{10}}\sum_ {j=1}^{n} \Big[c[J(dx_j)|\xi|^2+2\xi_{j}c[J(\xi)]\Big]
    \Big[\partial_{x_j}(c[J(\xi)])|\xi|^2-c[J(\xi)]\partial_{x_j}(|\xi|^2)\Big]\nonumber\\
    &+\frac{c[J(\xi)] T_J c[J(\xi)]}{|\xi|^6}-\frac{1}{|\xi|^8}\sum_ {j,l=1}^{n}c[J(\xi)]c[J(e_l)]c[J(e_j)] T_J c[J(\xi)] \xi_l \xi_j\nonumber\\
    &-\frac{1}{|\xi|^8}\sum_ {j,l=1}^{n}c[J(\xi)]c[J(e_l)]T_J c[J(e_j)] c[J(\xi)] \xi_l \xi_j\nonumber\\
    :=&Q_0+Q_1+Q_2\nonumber
\end{align}

\begin{align}
Q_0:=&\;\;\;\frac{c[J(\xi)]\sigma_{2}({D}_{J}^3)c[J(\xi)]}{|\xi|^8}\nonumber\\
&+\frac{c[J(\xi)]}{|\xi|^{10}}\sum_ {j=1}^{n} \Big[c[J(dx_j)|\xi|^2+2\xi_{j}c[J(\xi)]\Big]
\Big[\partial_{x_j}(c[J(\xi)])|\xi|^2-c[J(\xi)]\partial_{x_j}(|\xi|^2)\Big];\nonumber\\
Q_1:=&\frac{c[J(\xi)] T_J c[J(\xi)]}{|\xi|^6};\nonumber\\
Q_2:=&-\frac{1}{|\xi|^8}\sum_ {j,l=1}^{n}c[J(\xi)]c[J(e_l)]c[J(e_j)] T_J c[J(\xi)] \xi_l \xi_j-\frac{1}{|\xi|^8}\sum_ {j,l=1}^{n}c[J(\xi)]c[J(e_l)]T_J c[J(e_j)] c[J(\xi)] \xi_l \xi_j.\nonumber
\end{align}

A simple calculation gives
\begin{align}\label{Q0(x0)}
	&i\int_{|\xi'|=1}\int^{+\infty}_{-\infty}
    {\rm tr} [\partial_{\xi_n}\pi^+_{\xi_n}\sigma_{-1}({D}_{J,T}^{-1})\times Q_0](x_0)d\xi_n\sigma(\xi')dx'\\
    =&\frac{\pi}{16}\sum_{i=1}^{n-1}\left(a_{n}^{i}\right)^{2} \operatorname{tr}[\mathrm{id}] \Omega_{4} h^{\prime}(0) d x^{\prime}-\frac{\pi}{16}\sum_{i=1}^{n-1} a_{i}^{i} a_{n}^{n} \operatorname{tr}[\mathrm{id}] \Omega_{4} h^{\prime}(0) d x^{\prime} \nonumber\\
    &+\frac{11 \pi}{128}\sum_{l=1}^{n}\left(a_{n}^{l}\right)^{2} \operatorname{tr}[\mathrm{id}] \Omega_{4} h^{\prime}(0) d x^{\prime}-\frac{\pi}{64}\sum_{i=1}^{n-1}\left(a_{n}^{i}\right)^{2} \operatorname{tr}[\mathrm{id}] \Omega_{4} h^{\prime}(0) d x^{\prime}\nonumber\\
    &-\frac{23 \pi}{320}\sum_{l=1}^{n} \sum_{i=1}^{n-1}\left(a_{l}^{i}\right)^{2} \operatorname{tr}[\mathrm{id}] \Omega_{4} h^{\prime}(0)d x^{\prime}+\frac{\pi}{64}\sum_{i, j=1}^{n-1}\left(a_{j}^{i}\right)^{2} \operatorname{tr}[\mathrm{id}] \Omega_{4} h^{\prime}(0) d x^{\prime} \nonumber\\
    &+\frac{\pi}{64}\sum_{l, \alpha=1}^{n} \sum_{i=1}^{n-1}\left(a_{l}^{i}\right)^{2} g^{M}\left(J\left(d x_{\alpha}\right),\left(\nabla_{e_{\alpha}}^{L} J\right) e_{n}\right) \operatorname{tr}[\mathrm{id}] \Omega_{4} h^{\prime}(0) d x^{\prime} \nonumber\\
    &-\frac{\pi}{16}\sum_{i=1}^{n-1} g^{M}\left(J\left(d x_{i}\right),\left(\nabla_{e_{n}}^{L} J\right) e_{i}\right) \operatorname{tr}[\mathrm{id}] \Omega_{4} h^{\prime}(0)d x^{\prime} \nonumber\\
    &+\frac{\pi}{16}\sum_{i=1}^{n-1} g^{M}\left(J\left(d x_{n}\right),\left(\nabla_{e_{i}}^{L} J\right) e_{i}\right) \operatorname{tr}[\mathrm{id}] \Omega_{4} h^{\prime}(0) d x^{\prime} \nonumber\\
    &-\frac{\pi}{64}\sum_{l, \alpha=1}^{n}\left(a_{l}^{n}\right)^{2} g^{M}\left(J\left(d x_{\alpha}\right),\left(\nabla_{e_{\alpha}}^{L} J\right) e_{n}\right) \operatorname{tr}[\mathrm{id}] \Omega_{4} h^{\prime}(0) d x^{\prime} \nonumber\\
    &+\frac{\pi}{8}\sum_{l=1}^{n} \sum_{i=1}^{n-1} a_{l}^{i} \partial_{x_{i}}\left(a_{l}^{n}\right) \operatorname{tr}[\mathrm{id}] \Omega_{4} h^{\prime}(0) d x^{\prime}.\nonumber
\end{align}

\begin{align}\label{Q1}
	Q_1(x_0)|_{|{\xi}'|=1}=&\frac{c[J(\xi)]T_Jc[J(\xi)]}{|\xi|^6}(x_0)|_{|{\xi}'|=1}\\
	=&\frac{1}{4|\xi|^6} \bigg[\sum_{\alpha,q,v,l,r,u\neq s\neq t=1}^{n} \xi_{q}\xi_{\alpha}T_{vst} a_u^v a_l^q a_r^{\alpha} c(dx_l) c(e_u)c(e_s)c(e_t)c(dx_r) \nonumber\\
	&-2\sum_{\alpha,q,l,r,u,t,v=1}^{n} \xi_{q}\xi_{\alpha}T_{vut} a_u^v  a_l^q a_r^{\alpha} c(dx_l)c(e_t)c(dx_r)
	\bigg](x_0)|_{|{\xi}'|=1}\nonumber\\	
	=&\;\;\;\frac{\xi_n^2}{4(1+\xi_n^2)^3}\sum_{v,l,r,u\neq s\neq t=1}^{n}T_{vst} a_u^v a_l^n a_r^{n} c(dx_l) c(dx_u)c(dx_s)c(dx_t)c(dx_r)(x_0) \nonumber\\
	&+\frac{\xi_{n}}{4(1+\xi_n^2)^3}\sum_{v,l,r,u\neq s\neq t=1}^{n}\sum_{\alpha=1}^{n-1}\xi_{\alpha}T_{vst} a_u^v a_l^n a_r^{\alpha} c(dx_l) c(dx_u)c(dx_s)c(dx_t)c(dx_r)(x_0) \nonumber\\
	&+\frac{\xi_{n}}{4(1+\xi_n^2)^3}\sum_{v,l,r,u\neq s\neq t=1}^{n}\sum_{q=1}^{n-1}\xi_{q}T_{vst} a_u^v a_l^q a_r^{n} c(dx_l) c(dx_u)c(dx_s)c(dx_t)c(dx_r)(x_0) \nonumber\\
	&+\frac{1}{4(1+\xi_n^2)^3}\sum_{v,l,r,u\neq s\neq t=1}^{n}\sum_{\alpha,q=1}^{n-1}\xi_{q}\xi_{\alpha}T_{vst} a_u^v a_l^q a_r^{\alpha} c(dx_l) c(dx_u)c(dx_s)c(dx_t)c(dx_r)(x_0),\nonumber
\end{align}
and
\begin{align}\label{tr6}
 &\sum_{\beta,l,r,u\neq s\neq t=1}^{n}{\rm tr} \bigg[c(dx_{\beta}) c(dx_l) c(dx_u)c(dx_s)c(dx_t)c(dx_r)\bigg]\\
 =&\sum_{\beta,l,r,u\neq s\neq t=1}^{n}\bigg[-\delta_{\beta}^{t}\delta_{l}^{s}\delta_{u}^{r}+\delta_{\beta}^{t}\delta_{l}^{u}\delta_{s}^{r}+\delta_{\beta}^{s}\delta_{l}^{t}\delta_{u}^{r}-\delta_{\beta}^{s}\delta_{l}^{u}\delta_{t}^{r}-\delta_{\beta}^{u}\delta_{l}^{t}\delta_{s}^{r}+\delta_{\beta}^{u}\delta_{l}^{s}\delta_{t}^{r}\bigg].\nonumber
\end{align}
Obtained from \eqref{Q1} and \eqref{tr6}
\begin{align}\label{Q1(x0)}
	&i\int_{|\xi'|=1}\int^{+\infty}_{-\infty}
	{\rm tr} [\partial_{\xi_n}\pi^+_{\xi_n}\sigma_{-1}({D}_{J,T}^{-1})\times Q_1](x_0)d\xi_n\sigma(\xi')dx'\\
	=&i\int_{|\xi'|=1}\int^{+\infty}_{-\infty}-\frac{\xi_n}{8(\xi_n-i)(1+\xi_n^2)^2}\sum_{v,\beta,l,r,u\neq s\neq t=1}^{n}\sum_{i,\alpha=1}^{n-1}T_{vst}\xi_i \xi_{\alpha} a_{\beta}^i a_l^n a_r^{\alpha} a_u^v \nonumber\\
	&\times {\rm tr} \bigg[c(dx_{\beta}) c(dx_l) c(dx_u)c(dx_s)c(dx_t)c(dx_r) \bigg]\nonumber\\
	&-\frac{\xi_n}{8(\xi_n-i)(1+\xi_n^2)^2}\sum_{v,\beta,l,r,u\neq s\neq t=1}^{n}\sum_{i,q=1}^{n-1}T_{vst}\xi_i \xi_{q} a_{\beta}^i a_l^q a_r^{n} a_u^v \nonumber\\
	&\times{\rm tr} \bigg[c(dx_{\beta}) c(dx_l) c(dx_u)c(dx_s)c(dx_t)c(dx_r) \bigg]\nonumber\\
	&-\frac{i\xi_n^2}{8(\xi_n-i)(1+\xi_n^2)^2}\sum_{v,\beta,l,r,u\neq s\neq t=1}^{n}T_{vst} a_{\beta}^n a_l^n a_r^{n} a_u^v \nonumber\\
	&\times{\rm tr} \bigg[c(dx_{\beta}) c(dx_l) c(dx_u)c(dx_s)c(dx_t)c(dx_r) \bigg]\nonumber\\
	&-\frac{i}{8(\xi_n-i)(1+\xi_n^2)^2}\sum_{v,\beta,l,r,u\neq s\neq t=1}^{n}\sum_{\alpha,q=1}^{n-1}T_{vst} \xi_{\alpha}\xi_q a_{\beta}^n a_l^q a_r^{\alpha} a_u^v \nonumber\\
	&\times{\rm tr} \bigg[c(dx_{\beta}) c(dx_l) c(dx_u)c(dx_s)c(dx_t)c(dx_r) \bigg](x_0)d\xi_n\sigma(\xi')dx'\nonumber\\
	=&\;\;\;0.\nonumber
\end{align}

\begin{align}\label{Q2}
    Q_2(x_0)|_{|{\xi}'|=1}&=\bigg[-\frac{1}{|\xi|^8}\sum_ {j,l=1}^{n}c[J(\xi)]c[J(e_l)]c[J(e_j)] T_J c[J(\xi)] \xi_l \xi_j\\
    &\;\;\;-\frac{1}{|\xi|^8}\sum_ {j,l=1}^{n}c[J(\xi)]c[J(e_l)]T_J c[J(e_j)] c[J(\xi)] \xi_l \xi_j\bigg](x_0)|_{|{\xi}'|}\nonumber\\
    &=-\frac{1}{(1+\xi_n^2)^4}\sum_{q,\alpha,l,j,v,r,u\neq s\neq t=1}^{n}\xi_{q}\xi_{\alpha}\xi_{l}\xi_{j}T_{vst}a_p^q a_r^{\alpha} a_k^l a_h^j a_u^v \times \nonumber\\
    &\;\;\;\bigg[c(dx_p)c(dx_k)c(dx_h)c(dx_u)c(dx_s)c(dx_t)c(dx_r)\nonumber\\
    &\;\;\;+c(dx_p)c(dx_k)c(dx_u)c(dx_s)c(dx_t)c(dx_h)c(dx_r)\bigg](x_0),\nonumber
\end{align}

and

\begin{align}\label{tr8}
{\rm Tr}\aleph&:=\sum_{\beta,p,k,h,r,u\neq s\neq t=1}^{n}{\rm tr}\bigg[c(dx_\beta)c(dx_p)c(dx_k)c(dx_h)c(dx_u)c(dx_s)c(dx_t)c(dx_r)\\
&+c(dx_\beta)c(dx_p)c(dx_k)c(dx_u)c(dx_s)c(dx_t)c(dx_h)c(dx_r)\bigg]\nonumber\\
&=\sum_{\beta,p,k,h,r,u\neq s\neq t=1}^{n}\biggl\{ -2g(dx_h.dx_t){\rm tr}\bigg[c(dx_\beta)c(dx_p)c(dx_k)c(dx_u)c(dx_s)c(dx_r)\bigg]\nonumber\\
& \;\;\;+2g(dx_h.dx_s){\rm tr}\bigg[c(dx_\beta)c(dx_p)c(dx_k)c(dx_u)c(dx_t)c(dx_r)\bigg]\nonumber\\
&\;\;\;-2g(dx_h.dx_u){\rm tr}\bigg[c(dx_\beta)c(dx_p)c(dx_k)c(dx_s)c(dx_t)c(dx_r)\bigg]\biggr\}.\nonumber
\end{align}

Obtained from \eqref{Q2} and \eqref{tr8}
\begin{align}\label{Q2(x0)}
	&i\int_{|\xi'|=1}\int^{+\infty}_{-\infty}
	{\rm tr} [\partial_{\xi_n}\pi^+_{\xi_n}\sigma_{-1}({D}_{J,T}^{-1})\times Q_2](x_0)d\xi_n\sigma(\xi')dx'\\
	=&i\int_{|\xi'|=1}\int^{+\infty}_{-\infty}\biggl\{\frac{\xi_n^3}{8(\xi_n-i)^2(1+\xi_n^2)^4}\sum_{i,q=1}^{n-1}T_{vst}\xi_i \xi_{q} a_{\beta}^i a_p^q a_r^n a_k^n a_h^n a_u^v {\rm Tr}\aleph \nonumber\\
	&+\frac{\xi_n^3}{8(\xi_n-i)^2(1+\xi_n^2)^4}\sum_{i,\alpha=1}^{n-1}T_{vst}\xi_i \xi_{\alpha} a_{\beta}^i a_p^n a_r^{\alpha} a_k^n a_h^n a_u^v {\rm Tr}\aleph \nonumber\\
	&+\frac{\xi_n^3}{8(\xi_n-i)^2(1+\xi_n^2)^4}\sum_{i,l=1}^{n-1}T_{vst}\xi_i \xi_{l} a_{\beta}^i a_p^n a_r^{n} a_k^l a_h^n a_u^v {\rm Tr}\aleph \nonumber\\
	&+\frac{\xi_n^3}{8(\xi_n-i)^2(1+\xi_n^2)^4}\sum_{i,j=1}^{n-1}T_{vst}\xi_i \xi_{j} a_{\beta}^i a_p^n a_r^{n} a_k^n a_h^j a_u^v {\rm Tr}\aleph \nonumber\\
	&+\frac{\xi_n}{8(\xi_n-i)^2(1+\xi_n^2)^4}\sum_{i,q,\alpha,l=1}^{n-1}T_{vst}\xi_i \xi_q \xi_{\alpha} \xi_l a_{\beta}^i a_p^q a_r^{\alpha} a_k^l a_h^n a_u^v {\rm Tr}\aleph \nonumber\\	&+\frac{\xi_n}{8(\xi_n-i)^2(1+\xi_n^2)^4}\sum_{i,q,\alpha,j=1}^{n-1}T_{vst}\xi_i \xi_q \xi_{\alpha} \xi_j a_{\beta}^i a_p^q a_r^{\alpha} a_k^n a_h^j a_u^v {\rm Tr}\aleph \nonumber\\
	&+\frac{\xi_n}{8(\xi_n-i)^2(1+\xi_n^2)^4}\sum_{i,q,l,j=1}^{n-1}T_{vst}\xi_i \xi_q \xi_{l} \xi_j a_{\beta}^i a_p^q a_r^{n} a_k^l a_h^j a_u^v {\rm Tr}\aleph \nonumber\\		&+\frac{\xi_n}{8(\xi_n-i)^2(1+\xi_n^2)^4}\sum_{i,\alpha,l,j=1}^{n-1}T_{vst}\xi_i \xi_\alpha \xi_{l} \xi_j a_{\beta}^i a_p^n a_r^{\alpha} a_k^l a_h^j a_u^v {\rm Tr}\aleph \nonumber\\
	&+\frac{i\xi_n^4}{8(\xi_n-i)^2(1+\xi_n^2)^4}\sum_{i,\alpha,l,j=1}^{n-1}T_{vst} a_{\beta}^n a_p^n a_r^{n} a_k^n a_h^n a_u^v {\rm Tr}\aleph \nonumber\\
	&+\frac{i\xi_n^2}{8(\xi_n-i)^2(1+\xi_n^2)^4}\sum_{q,\alpha=1}^{n-1}T_{vst}\xi_q \xi_\alpha a_{\beta}^n a_p^q a_r^{\alpha} a_k^n a_h^n a_u^v {\rm Tr}\aleph \nonumber\\
	&+\frac{i\xi_n^2}{8(\xi_n-i)^2(1+\xi_n^2)^4}\sum_{q,l=1}^{n-1}T_{vst}\xi_q \xi_l a_{\beta}^n a_p^q a_r^{n} a_k^l a_h^n a_u^v {\rm Tr}\aleph \nonumber\\
	&+\frac{i\xi_n^2}{8(\xi_n-i)^2(1+\xi_n^2)^4}\sum_{q,j=1}^{n-1}T_{vst}\xi_q \xi_j a_{\beta}^n a_p^q a_r^{n} a_k^n a_h^j a_u^v {\rm Tr}\aleph \nonumber\\
	&+\frac{i\xi_n^2}{8(\xi_n-i)^2(1+\xi_n^2)^4}\sum_{\alpha,l=1}^{n-1}T_{vst}\xi_\alpha \xi_l a_{\beta}^n a_p^n a_r^{\alpha} a_k^l a_h^n a_u^v {\rm Tr}\aleph \nonumber\\
	&+\frac{i\xi_n^2}{8(\xi_n-i)^2(1+\xi_n^2)^4}\sum_{\alpha,j=1}^{n-1}T_{vst}\xi_\alpha \xi_j a_{\beta}^n a_p^n a_r^{\alpha} a_k^n a_h^j a_u^v {\rm Tr}\aleph \nonumber\\
	&+\frac{i\xi_n^2}{8(\xi_n-i)^2(1+\xi_n^2)^4}\sum_{l,j=1}^{n-1}T_{vst}\xi_l \xi_j a_{\beta}^n a_p^n a_r^{n} a_k^l a_h^j a_u^v {\rm Tr}\aleph \biggr\}(x_0)d\xi_n\sigma(\xi')dx'\nonumber\\
	=&\;\;\;0.\nonumber
\end{align}

Combining \eqref{Q0(x0)} and \eqref{Q1(x0)} and \eqref{Q2(x0)} yields
\begin{align}
	\hat{\Psi_4}&=\frac{\pi}{16}\sum_{i=1}^{n-1}\left(a_{n}^{i}\right)^{2} \operatorname{tr}[\mathrm{id}] \Omega_{4} h^{\prime}(0) d x^{\prime}-\frac{\pi}{16}\sum_{i=1}^{n-1} a_{i}^{i} a_{n}^{n} \operatorname{tr}[\mathrm{id}] \Omega_{4} h^{\prime}(0) d x^{\prime} \nonumber\\
	&+\frac{11 \pi}{128}\sum_{l=1}^{n}\left(a_{n}^{l}\right)^{2} \operatorname{tr}[\mathrm{id}] \Omega_{4} h^{\prime}(0) d x^{\prime}-\frac{\pi}{64}\sum_{i=1}^{n-1}\left(a_{n}^{i}\right)^{2} \operatorname{tr}[\mathrm{id}] \Omega_{4} h^{\prime}(0) d x^{\prime}\nonumber\\
	&-\frac{23 \pi}{320}\sum_{l=1}^{n} \sum_{i=1}^{n-1}\left(a_{l}^{i}\right)^{2} \operatorname{tr}[\mathrm{id}] \Omega_{4} h^{\prime}(0)d x^{\prime}+\frac{\pi}{64}\sum_{i, j=1}^{n-1}\left(a_{j}^{i}\right)^{2} \operatorname{tr}[\mathrm{id}] \Omega_{4} h^{\prime}(0) d x^{\prime} \nonumber\\
	&+\frac{\pi}{64}\sum_{l, \alpha=1}^{n} \sum_{i=1}^{n-1}\left(a_{l}^{i}\right)^{2} g^{M}\left(J\left(d x_{\alpha}\right),\left(\nabla_{e_{\alpha}}^{L} J\right) e_{n}\right) \operatorname{tr}[\mathrm{id}] \Omega_{4} h^{\prime}(0) d x^{\prime} \nonumber\\
	&-\frac{\pi}{16}\sum_{i=1}^{n-1} g^{M}\left(J\left(d x_{i}\right),\left(\nabla_{e_{n}}^{L} J\right) e_{i}\right) \operatorname{tr}[\mathrm{id}] \Omega_{4} h^{\prime}(0)d x^{\prime} \nonumber\\
	&+\frac{\pi}{16}\sum_{i=1}^{n-1} g^{M}\left(J\left(d x_{n}\right),\left(\nabla_{e_{i}}^{L} J\right) e_{i}\right) \operatorname{tr}[\mathrm{id}] \Omega_{4} h^{\prime}(0) d x^{\prime} \nonumber\\
	&-\frac{\pi}{64}\sum_{l, \alpha=1}^{n}\left(a_{l}^{n}\right)^{2} g^{M}\left(J\left(d x_{\alpha}\right),\left(\nabla_{e_{\alpha}}^{L} J\right) e_{n}\right) \operatorname{tr}[\mathrm{id}] \Omega_{4} h^{\prime}(0) d x^{\prime} \nonumber\\
	&+\frac{\pi}{8}\sum_{l=1}^{n} \sum_{i=1}^{n-1} a_{l}^{i} \partial_{x_{i}}\left(a_{l}^{n}\right) \operatorname{tr}[\mathrm{id}] \Omega_{4} h^{\prime}(0) d x^{\prime}.\nonumber
\end{align}

\noindent  { \bf case $\bf \hat{c}$)}~$r=-2,~l=-3,~k=j=|\alpha|=0$\\
\noindent By (3.18), we calculate that
\begin{align}
	\hat{\Psi_5}=&-i\int_{|\xi'|=1}\int^{+\infty}_{-\infty}
	{\rm tr} [\pi^+_{\xi_n}\sigma_{-2}({D}_{J,T}^{-1})\times
	\partial_{\xi_n}\sigma_{-3}({D}_{J,T}^{-3})](x_0)d\xi_n\sigma(\xi')dx'\nonumber
\end{align}

Same classification as in case b \eqref{B}
By computation, we have
\begin{align}\label{6B_0}
	&-i\int_{|\xi'|=1}\int^{+\infty}_{-\infty}{\rm tr} [\pi^+_{\xi_n}B_0(x_0)\times
    \partial_{\xi_n}\sigma_{-3}({D}_{J,T}^{-3})](x_0)d\xi_n\sigma(\xi')dx'\\
    & =-\frac{\pi}{64}\sum_{l, j, \beta=1}^{n}\left(a_{\beta}^{n}\right)^{2} a_{l}^{j} \partial_{x_{j}}\left(a_{l}^{n}\right) \operatorname{tr}[\mathrm{id}] \Omega_{4} d x^{\prime}-\frac{\pi}{64}\sum_{l, j, \beta=1}^{n} \sum_{i=1}^{n-1}\left(a_{\beta}^{i}\right)^{2} a_{l}^{j} \partial_{x_{j}}\left(a_{l}^{n}\right) \operatorname{tr}[\mathrm{id}] \Omega_{4}d x^{\prime} \nonumber\\
    &\;\;\; +\frac{\pi}{256}\sum_{l=1}^{n} \sum_{i=1}^{n-1}\left(a_{l}^{n}\right)^{2} a_{i}^{i} a_{n}^{n} \operatorname{tr}[\mathrm{id}] \Omega_{4} h^{\prime}(0) d x^{\prime}+\frac{3 \pi}{256}\sum_{l=1}^{n} \sum_{\nu, i=1}^{n-1}\left(a_{l}^{i}\right)^{2} a_{\nu}^{\nu} a_{n}^{n} \operatorname{tr}[\mathrm{id}] \Omega_{4} h^{\prime}(0) d x^{\prime} \nonumber\\
    & \;\;\;-\frac{\pi}{32}\sum_{l, j=1}^{n} a_{l}^{j} \partial_{x_{j}}\left(a_{l}^{n}\right) \operatorname{tr}[\mathrm{id}] \Omega_{4} d x^{\prime}+\frac{\pi}{256}\sum_{l, \beta=1}^{n}\left(a_{\beta}^{n}\right)^{2}\left(a_{l}^{n}\right)^{2} \operatorname{tr}[\mathrm{id}] \Omega_{4} h^{\prime}(0)d x^{\prime} \nonumber\\
    &\;\;\; +\frac{49 \pi}{1280}\sum_{l, \beta=1}^{n} \sum_{i=1}^{n-1}\left(a_{l}^{n}\right)^{2}\left(a_{\beta}^{i}\right)^{2} \operatorname{tr}[\mathrm{id}] \Omega_{4} h^{\prime}(0)d x^{\prime}-\frac{3 \pi}{256}\sum_{l=1}^{n} \sum_{i=1}^{n-1}\left(a_{i}^{n}\right)^{2}\left(a_{l}^{n}\right)^{2} \operatorname{tr}[\mathrm{id}] \Omega_{4} h^{\prime}(0) d x^{\prime} \nonumber\\
    &\;\;\; -\frac{5 \pi}{256}\sum_{l=1}^{n} \sum_{\nu, i=1}^{n-1}\left(a_{\nu}^{n}\right)^{2}\left(a_{l}^{i}\right)^{2} \operatorname{tr}[\mathrm{id}] \Omega_{4} h^{\prime}(0) d x^{\prime}-\frac{\pi}{64}\sum_{l=1}^{n} \sum_{\nu, i=1}^{n-1}\left(a_{l}^{n}\right)^{2}\left(a_{\nu}^{i}\right)^{2} \operatorname{tr}[\mathrm{id}] \Omega_{4} h^{\prime}(0) d x^{\prime}\nonumber
\end{align}
According to \eqref{tr6} and \eqref{B} , the following equation holds
\begin{align}\label{6B_1}
	&-i\int_{|\xi'|=1}\int^{+\infty}_{-\infty}{\rm tr} [\pi^+_{\xi_n}B_1(x_0)\times
	\partial_{\xi_n}\sigma_{-3}({D}_{J,T}^{-3})](x_0)d\xi_n\sigma(\xi')dx'\\
    =&-i\int_{|\xi'|=1}\int^{+\infty}_{-\infty}\nonumber\\
    &\biggl\{-\frac{\xi_{n}}{4(\xi_{n}-i)^2(1+\xi_{n}^2)^3 }\sum_{\beta, l,r,v,u\neq s\neq t=1}^{n} \sum_{i,\alpha=1}^{n-1} T_{vst}\xi_{i}\xi_{\alpha} a_l^n a_r^{\alpha} a_u^v  a_{\beta}^{i} {\rm tr} [c(dx_l) c(dx_u) c(dx_s) c(dx_t) c(dx_r)c(dx_{\beta})]\nonumber\\
    &-\frac{\xi_{n}}{4(\xi_{n}-i)^2(1+\xi_{n}^2)^3} \sum_{\beta, l,r,v,u\neq s\neq t=1}^{n}\sum_{i,q=1}^{n-1} T_{vst} \xi_{i}\xi_{q} a_l^q a_r^n a_u^v  a_{\beta}^{i} {\rm tr} [c(dx_l) c(dx_u) c(dx_s) c(dx_t)  c(dx_r)c(dx_{\beta})]\nonumber\\
    &+\frac{\xi_{n}-3\xi_{n}^{3}}{16(\xi_{n}-i)^2(1+\xi_{n}^2)^3}\sum_{\beta, l,r,v,u\neq s\neq t=1}^{n}  T_{vst}  a_l^n a_r^n a_u^v a_{\beta}^{n} {\rm tr} [c(dx_l) c(dx_u) c(dx_s) c(dx_t) c(dx_r)c(dx_{\beta})]\nonumber\\
    &-\frac{i(1-3\xi_{n}^2)(i \xi_{n}+2)}{16(\xi_{n}-i)^2(1+\xi_{n}^2)^3 } \sum_{\beta, l,r,v,u\neq s\neq t=1}^{n}\sum_{q,\alpha=1}^{n-1} T_{vst} \xi_{q} \xi_{\alpha} a_l^q a_r^{\alpha} a_u^v  a_{\beta}^{n} {\rm tr} [c(dx_l) c(dx_u) c(dx_s) c(dx_t) c(dx_r)c(dx_{\beta})]\biggr\}\nonumber\\
    &(x_0)d\xi_n\sigma(\xi')dx'\nonumber\\
    &=0.\nonumber
\end{align}
\begin{align}\label{6B_2}
	&-i\int_{|\xi'|=1}\int^{+\infty}_{-\infty}{\rm tr} [\pi^+_{\xi_n}B_2(x_0)\times
	\partial_{\xi_n}\sigma_{-3}({D}_{J,T}^{-3})](x_0)d\xi_n\sigma(\xi')dx'\\
	=&-i\int_{|\xi'|=1}\int^{+\infty}_{-\infty}\biggl\{\frac{\xi_{n}}{2(\xi_{n}-i)^2(1+\xi_{n}^2)^3 }\sum_{\beta, l,r,v,u,t=1}^{n} \sum_{i,\alpha=1}^{n-1} T_{vut}\xi_{i}\xi_{\alpha} a_l^n a_r^{\alpha} a_u^v  a_{\beta}^{i} {\rm tr} [c(dx_l)c(dx_t) c(dx_r)c(dx_{\beta})]\nonumber\\
	&+\frac{\xi_{n}}{2(\xi_{n}-i)^2(1+\xi_{n}^2)^3} \sum_{\beta, l,r,v,u,t=1}^{n}\sum_{i,q=1}^{n-1} T_{vut} \xi_{i}\xi_{q} a_l^q a_r^n a_u^v  a_{\beta}^{i} {\rm tr} [c(dx_l) c(dx_t)  c(dx_r)c(dx_{\beta})]\nonumber\\
	&+\frac{i(1-3\xi_{n}^2)(i \xi_{n}+2)}{8(\xi_{n}-i)^2(1+\xi_{n}^2)^3 } \sum_{\beta, l,r,v,u,t=1}^{n}\sum_{q,\alpha=1}^{n-1} T_{vut} \xi_{q} \xi_{\alpha} a_l^q a_r^{\alpha} a_u^v  a_{\beta}^{n} {\rm tr} [c(dx_l) c(dx_t) c(dx_r)c(dx_{\beta})]\nonumber\\	&-\frac{\xi_{n}-3\xi_{n}^{3}}{8(\xi_{n}-i)^2(1+\xi_{n}^2)^3}\sum_{\beta, l,r,v,u,t=1}^{n}  T_{vut}  a_l^n a_r^n a_u^v a_{\beta}^{n} {\rm tr} [c(dx_l) c(dx_t) c(dx_r)c(dx_{\beta})]\biggr\}(x_0)d\xi_n\sigma(\xi')dx'\nonumber\\
	=&-i\int_{|\xi'|=1}\int^{+\infty}_{-\infty}\biggl\{\frac{\xi_{n}}{2(\xi_{n}-i)^2(1+\xi_{n}^2)^3 }\sum_{\beta, l,r,v,u,t=1}^{n} \sum_{i,\alpha=1}^{n-1} T_{vut}\xi_{i}\xi_{\alpha} a_t^n a_r^{\alpha} a_u^v  a_{r}^{i} \nonumber\\
	&+\frac{\xi_{n}}{2(\xi_{n}-i)^2(1+\xi_{n}^2)^3} \sum_{\beta, l,r,v,u,t=1}^{n}\sum_{i,q=1}^{n-1} T_{vut} \xi_{i}\xi_{q} a_l^q a_t^n a_u^v  a_{l}^{i}\nonumber\\
	&-\frac{i(1-3\xi_{n}^2)(i \xi_{n}+2)}{8(\xi_{n}-i)^2(1+\xi_{n}^2)^3 } \sum_{\beta, l,r,v,u,t=1}^{n}\sum_{q,\alpha=1}^{n-1} T_{vut} \xi_{q} \xi_{\alpha} a_l^q a_l^{\alpha} a_u^v  a_{t}^{n} \nonumber\\	
	&-\frac{\xi_{n}-3\xi_{n}^{3}}{8(\xi_{n}-i)^2(1+\xi_{n}^2)^3}\sum_{\beta, l,r,v,u,t=1}^{n}  T_{vut}  a_t^n a_r^n a_u^v a_{r}^{n} \biggr\}(x_0)d\xi_n\sigma(\xi')dx'\nonumber\\
	&=0.\nonumber
\end{align}

\eqref{6B_0}-\eqref{6B_2}
\begin{align}
	\hat{\Psi_5}& =-\frac{\pi}{64}\sum_{l, j, \beta=1}^{n}\left(a_{\beta}^{n}\right)^{2} a_{l}^{j} \partial_{x_{j}}\left(a_{l}^{n}\right) \operatorname{tr}[\mathrm{id}] \Omega_{4} d x^{\prime}-\frac{\pi}{64}\sum_{l, j, \beta=1}^{n} \sum_{i=1}^{n-1}\left(a_{\beta}^{i}\right)^{2} a_{l}^{j} \partial_{x_{j}}\left(a_{l}^{n}\right) \operatorname{tr}[\mathrm{id}] \Omega_{4}d x^{\prime} \nonumber\\
    &\;\;\; +\frac{\pi}{256}\sum_{l=1}^{n} \sum_{i=1}^{n-1}\left(a_{l}^{n}\right)^{2} a_{i}^{i} a_{n}^{n} \operatorname{tr}[\mathrm{id}] \Omega_{4} h^{\prime}(0) d x^{\prime}+\frac{3 \pi}{256}\sum_{l=1}^{n} \sum_{\nu, i=1}^{n-1}\left(a_{l}^{i}\right)^{2} a_{\nu}^{\nu} a_{n}^{n} \operatorname{tr}[\mathrm{id}] \Omega_{4} h^{\prime}(0) d x^{\prime} \nonumber\\
    & \;\;\;-\frac{\pi}{32}\sum_{l, j=1}^{n} a_{l}^{j} \partial_{x_{j}}\left(a_{l}^{n}\right) \operatorname{tr}[\mathrm{id}] \Omega_{4} d x^{\prime}+\frac{\pi}{256}\sum_{l, \beta=1}^{n}\left(a_{\beta}^{n}\right)^{2}\left(a_{l}^{n}\right)^{2} \operatorname{tr}[\mathrm{id}] \Omega_{4} h^{\prime}(0)d x^{\prime} \nonumber\\
    &\;\;\; +\frac{49 \pi}{1280}\sum_{l, \beta=1}^{n} \sum_{i=1}^{n-1}\left(a_{l}^{n}\right)^{2}\left(a_{\beta}^{i}\right)^{2} \operatorname{tr}[\mathrm{id}] \Omega_{4} h^{\prime}(0)d x^{\prime}-\frac{3 \pi}{256}\sum_{l=1}^{n} \sum_{i=1}^{n-1}\left(a_{i}^{n}\right)^{2}\left(a_{l}^{n}\right)^{2} \operatorname{tr}[\mathrm{id}] \Omega_{4} h^{\prime}(0) d x^{\prime} \nonumber\\
    &\;\;\; -\frac{5 \pi}{256}\sum_{l=1}^{n} \sum_{\nu, i=1}^{n-1}\left(a_{\nu}^{n}\right)^{2}\left(a_{l}^{i}\right)^{2} \operatorname{tr}[\mathrm{id}] \Omega_{4} h^{\prime}(0) d x^{\prime}-\frac{\pi}{64}\sum_{l=1}^{n} \sum_{\nu, i=1}^{n-1}\left(a_{l}^{n}\right)^{2}\left(a_{\nu}^{i}\right)^{2} \operatorname{tr}[\mathrm{id}] \Omega_{4} h^{\prime}(0) d x^{\prime}.\nonumber
\end{align}

By simple calculation
\begin{align}
\sum_{i=1}^{n-1} g^{M}\left(J\left(d x_{i}\right),\left(\nabla_{e_{n}}^{L} J\right) e_{i}\right)\left(x_{0}\right)=\sum_{i=1}^{n-1} g^{M}\left(J\left(d x_{i}\right), \nabla_{e_{n}}^{L} J\left(e_{i}\right)-J\left(\nabla_{e_{n}}^{L} e_{i}\right)\right)\left(x_{0}\right)=0,\nonumber
\end{align}
by similar calculations, we can get
\begin{align}
\sum_{i=1}^{n-1} g^{M}\left(J\left(d x_{n}\right),\left(\nabla_{e_{i}}^{L} J\right) e_{i}\right)\left(x_{0}\right)=0 ;\nonumber\\
\sum_{\alpha=1}^{n} g^{M}\left(J\left(d x_{\alpha}\right),\left(\nabla_{e_{\alpha}}^{L} J\right) e_{n}\right)\left(x_{0}\right)=0.\nonumber
\end{align}

In summary,
\begin{align}
	\Psi=\Psi_1+\Psi_2+\Psi_3+\Psi_4+\Psi_5=-\frac{\pi}{2}\big(1-\left\langle J(e_6),e_6 \right\rangle ^2\big)h^{\prime}(0)\Omega_{4}  d x^{\prime}.\nonumber
\end{align}

\begin{thm}
	Let $M$ be a $6$-dimensional almost product Riemannian spin manifold with the boundary $\partial M$ and the metric
	$g^M$ as above, ${{D}_{J,T}}$ be the $J$-twist of the Dirac operator with torsion on $\widetilde{M}$, then
	\begin{align}
		&\widetilde{{\rm Wres}}[\pi^+{{D}_{J,T}}^{-1}\circ\pi^+{{D}_{J,T}}^{-1}]\\
		=&\int_{M}4\pi^{2}\biggl\{ \sum_{i,j=1}^{n}R(J(e_{i}), J(e_{j}), e_{j}, e_{i})
		-2\sum_{\nu,j=1}^{n}g^{M}(\nabla_{e_{j}}^{L}(J)e_{\nu}, (\nabla^{L}_{e_{\nu}}J)e_{j})\nonumber \\
		&-2\sum_{\nu,j=1}^{n}g^{M}(J(e_{\nu}), (\nabla^{L}_{e_{j}}(\nabla^{L}_{e_{\nu}}(J)))e_{j}-(\nabla^{L}_{\nabla^{L}_{e_{j}}e_{\nu}}(J))e_{j})\nonumber\\
		&-\sum_{\alpha,\nu,j=1}^{n}g^{M}(J(e_{\alpha}), (\nabla^{L}_{e_{\nu}}J)e_{j})g^{M}((\nabla^{L}_{e_{\alpha}}J)e_{j}, J(e_{\nu}))\nonumber\\
	&-\sum_{\alpha,\nu,j=1}^{n}g^{M}(J(e_{\alpha}), (\nabla^{L}_{e_{\alpha}}J)e_{j})g^{M}(J(e_{\nu}), (\nabla^{L}_{e_{\nu}}J)e_{j})+\sum_{\nu,j=1}^{n}g^{M}((\nabla^{L}_{e_{\nu}}J)e_{j}, (\nabla^{L}_{e_{\nu}}J)e_{j})-\frac{5}{3}s\nonumber\\
	&-2\sum_{j, m, p} T\left(e_m, e_{p}, J\left(e_{m}\right)\right) g^{M}\left((\nabla_{e_{j}}^{L}J) e_{j} ,e_{p}\right)-2\sum_{j, m, \alpha} T\left(e_m, J\left(e_{j}\right), J\left(e_{\alpha}\right)\right) g^{M}\left((\nabla_{e_{\alpha}}^{L}J) e_{j} ,J(e_m)\right)\nonumber\\
	&-2\sum_{j, m, \alpha} T\left(e_m, J\left(e_{j}\right), J\left(e_{m}\right)\right) g^{M}\left((\nabla_{e_{\alpha}}^{L}J) e_{j} ,J(e_{\alpha})\right)+2\sum_{m,p, \alpha} T\left(e_m, e_p, J\left(e_{\alpha}\right)\right) g^{M}\left((\nabla_{e_{\alpha}}^{L}J) e_{m} ,e_p\right)\nonumber\\
	&-2\sum_{m,p, \alpha} T\left(e_m, e_p, J\left(e_{\alpha}\right)\right) g^{M}\left((\nabla_{e_{m}}^{L}J) e_{\alpha} ,e_p\right)- \sum_{j,p,m} T(e_j, e_p, J(e_m))T(e_m, e_p, J(e_j))\nonumber\\
	&-2 \sum_{j, l, \alpha} T^2\left(J\left(e_{j}\right), e_{l}, e_{\alpha}\right)\biggr\}d{\rm Vol_M}-\frac{\pi}{2}\int_{\partial M}\left(1- \left\langle J\left(e_{6}\right), e_{6}\right\rangle^{2}\right) h^{\prime}(0)\Omega_{4} d \mathrm{Vol}_{\partial \mathrm{M}}\nonumber
	\end{align}
	where $s$ is the scalar curvature.
\end{thm}

\vskip 1 true cm

\section{Declarations}
Ethics approval and consent to participate No applicable.\\

Consent for publication No applicable.\\

Availability of data and material The authors confirm that the data supporting the findings of this study are available within the article.\\

Competing interests The authors declare no conflict of interest.\\

Funding Sponsored by Natural Science Foundation of Xinjiang Uygur AutonomousRegion 2024D01C341 and the National Natural Science Foundation of China 11771070, 12061078.\\

Authors' contributions All authors contributed to the study conception and design. Material preparation, data collection and analysis were performed by Jin Hong, Siyao Liu and Yong Wang. The first draft of the manuscript was written by Jin Hong and all authors commented on previous versions of the manuscript. All authors read and approved the final manuscript.\\

\vskip 1 true cm


\bigskip
\bigskip

\noindent {\footnotesize {\it J. Hong} \\
	{School of Mathematics and Statistics, Yili Normal University, Yining 835000, China}\\
	{School of Mathematics and Statistics, Northeast Normal University, Changchun 130024, China}\\
	{Email: jhong@nenu.edu.cn}

\noindent {\footnotesize {\it S. Liu} \\
{School of Mathematics and Statistics, Changchun University of Science and Technology, Changchun
	130022, China}\\
{Email: liusy719@nenu.edu.cn}

\noindent {\footnotesize {\it Y. Wang} \\
{School of Mathematics and Statistics, Northeast Normal University, Changchun 130024, China}\\
{Email: wangy581@nenu.edu.cn}

\end{document}